\documentclass[11pt]{article}
\usepackage{amssymb}
\usepackage{amsmath}
\usepackage{amsthm}
\usepackage{float}
\usepackage{hyperref}

\newtheoremstyle{danielstheoremstyle}
  {\topsep} 
  {\topsep} 
  {\it} 
  {15pt} 
  {\bfseries} 
  {.} 
  {.5em} 
  {} 
\newtheoremstyle{danielsclaimstyle}
  {\topsep} 
  {\topsep} 
  {\it} 
  {} 
  {\bfseries} 
  {.} 
  {.5em} 
  {} 
\newtheoremstyle{danielsdefinitionstyle}
  {\topsep} 
  {\topsep} 
  {} 
  {15pt} 
  {\bfseries} 
  {.} 
  {.5em} 
  {} 
\newtheoremstyle{danielsremarkstyle}
  {\topsep} 
  {\topsep} 
  {} 
  {} 
  {\it} 
  {.} 
  {.5em} 
  {} 

\theoremstyle{danielstheoremstyle}
\newtheorem{theorem}{Theorem}[section]
\newtheorem{corollary}[theorem]{Corollary}
\newtheorem{lemma}[theorem]{Lemma}
\newtheorem{proposition}[theorem]{Proposition}
\theoremstyle{danielsclaimstyle}
\newtheorem{claim}[theorem]{Claim}
\theoremstyle{danielsdefinitionstyle}
\newtheorem{definition}[theorem]{Definition}
\theoremstyle{danielsremarkstyle}
\newtheorem{remark}[theorem]{Remark}
\newtheorem{example}[theorem]{Example}

\newenvironment{myremark}{\vspace{-25pt} \begin{quote} \begin{remark}}{\end{remark} \end{quote}}

\DeclareMathOperator\Log{Log}
\newcommand{\tngvec}[3][\gamma]{\frac{\partial}{\partial #1_{#2}}\Big\rvert_{#3}}
\newcommand{\quessianobb}[2][\frac{d'}{2}]{ \mathcal{Q}_{#1}\big(#2\big) \left( \frac{a(\xi)}{(1-|\xi|^{2})^{n+1}} \right)}
\newcommand{\detr}[1]{ \mathrm{det}_{\mathcal{R}} \mathcal{O}_{#1} }

\newcommand{\norma}[2]{\left\vert \left\vert #1 \right\vert \right\vert_{\scriptscriptstyle #2}}
\newcommand{\topnorm}{\widehat{T}_{a\mathrm{d}\sigma}}
\newcommand{\cirquess}[1][\varphi]{ \mathring{\mathcal{Q}}\left( #1 \right) }
\newcommand{\sdots}{ {\scriptstyle \dots } }
\newcommand{\scdots}{ {\scriptstyle \cdots } }

\newcommand{\normaaa}[2]{\big\vert \big\vert #1 \big\vert \big\vert_{\scriptscriptstyle #2}}
\allowdisplaybreaks

\begin{document}

\begin{center}
\textbf{ 
A Szeg\"{o} limit theorem for a class of Toeplitz operators on the Bergman space of the unit ball with singular symbols
}\\
\begin{small}
Daniel Ivan Ramirez Monta\~no\footnote{Instituto de Matem\'{a}ticas, Unidad Cuernavaca, Universidad Nacional Aut\'{o}noma de M\'{e}xico. {\it E-mail address}: daniel.ramirez@im.unam.mx}
\end{small}
\end{center}

\renewcommand{\abstractname}{\vspace{-\baselineskip}}
\begin{abstract}
We obtain a Szeg\"{o} limit theorem for a family of Toeplitz operators defined on the weighted Bergman space of the unit ball $\mathbb{B}_{n}$. The symbols of these operators are supported on some isotropic or co-isotropic submanifold $\Gamma \subseteq \mathbb{B}_{n}$ and can be seen, in general, as measures that are singular with respect to the Lebesgue measure on $\mathbb{C}^{n}$. The given theorem allows to describe the asymptotic behavior of these operators as the parameter of the weighted Bergman space tends to infinity.
\end{abstract}
\tableofcontents

\section{Introduction}\label{laintroduccion}

In general, a Szeg\"{o} limit theorem is a result that describes the asymptotic behavior of the spectral measures of an indexed family of operators as the parameter that indexes the family tends to a limit. This description is often made in terms of operator traces. Originally, G. Szeg\"{o} delivered expressions for the limit of determinants of truncations of infinite Toeplitz matrices as the truncation order tended to infinity (\cite{szego1}). These relations were later used to obtain information on the asymptotic distribution of the eigenvalues of the said matrices, expressing the limit of their traces through integrals. In particular, in \cite{grenander} it is shown that 
\begin{equation}\label{eqenintrodos}
\lim\limits_{N \rightarrow \infty} \frac{1}{N} \sum_{j=1}^{N} f(\lambda_{j,N}) = \frac{1}{2\pi} \int_{0}^{2\pi} f(a(x)) \ \mathrm{d}x \ , \nonumber
\end{equation}
where $\lambda_{j,N}$ are the eigenvalues of the $N$-th order truncation of the Toeplitz matrix whose diagonals are the Fourier coefficients of a real valued $ a \in L^{\infty}(\mathbb{S}^{1}) $ and $f $ is any continuous function on $[\mathrm{ess} \inf a , \mathrm{ess} \sup a]$. Later works focused on reproducing this result with weaker conditions for $a$ and applying to different families of functions $f$ (e.g. \cite{widom}, \cite{rusos}).\par

The generalization of some of these ideas to the operator setting is described in \cite{guillemin}. For a bounded self-adjoint operator $A$ on a Hilbert space $H$, we want a measure $\mu_{A}$ on $\mathbb{R}$ such that \begin{equation}\label{eqenintroigualdadtrazas}
\mathrm{Tr}( f(A) ) = \int_{\mathbb{R}} f \ \mathrm{d}\mu_{A} ,
\end{equation}
for any $f \in C(\mathbb{R})$. If $H$ has finite dimension, this can be done using the Riesz-Markov-Kakutani representation theorem, since the assignation $f \mapsto \mathrm{Tr}(f(A))$ results in a continuous functional on $C(\mathbb{R})$. In fact, $\mu_{A}$ is the sum of the Dirac measures centered on each eigenvalue of $A$. In general, this cannot be done when $H$ is not of finite dimension, so an alternative construction is proposed. Taking an increasing sequence of finite-rank orthogonal projections $P_{n}$ on $H$ such that $ P_{n}$ converges to the identity in $H$, the correspondent sequence of measures $\mu_{n}$ satisfying (\ref{eqenintroigualdadtrazas}) for the operator $P_{n}A\vert_{P_{n}(H)} $ is constructed. Therefore, a natural substitute for $\mu_{A}$ would be the weak limit of the $\mu_{n}$ measures. In this situation a Szeg\"{o} limit theorem translates to finding a measure whose integral reproduces the limit of the traces of continuous images of a sequence of operators. Particular attention has been paid to this problem when $A$ is replaced by a pseudo-differential operator (see \cite{guillemin}, \cite{laptev}, \cite{widomtres}). \par
An interesting scenario occurs when considering a weighted Hilbert space of complex functions as the dominion of the Toeplitz operators. By varying the parameter that defines the associated weights we obtain a parametrized family of operators whose limiting behavior can be studied. This idea was followed in both \cite{elnuevo} and \cite{SalvadorAlejandroUribe}. In \cite{elnuevo} results of this type were obtained considering Toeplitz operators on the Bergman space of the disc $\mathbb{D}$ and the Fock space of $\mathbb{C}$. Previously in \cite{SalvadorAlejandroUribe} the problem was studied in the context of the Fock space of $\mathbb{C}^{n}$ and considering symbols supported on some submanifold of $\mathbb{C}^{n}$, thus in general being measures that are singular with respect to the Lebesgue measure. We will also follow this approach, obtaining an analogue result to the main one of \cite{SalvadorAlejandroUribe} for the case of the weighted Bergman space of the unit ball $\mathbb{B}_{n}$.


\subsection{Definitions}

For $n \in \mathbb{N}$, let $\mathbb{B}_{n}:= \{ z \in \mathbb{C}^{n} : |z|<1 \}$ denote the unit ball in $\mathbb{C}^{n}$. For $\alpha>-1$, the weighted Bergman space $A^{2}_{\alpha}(\mathbb{B}_{n})$ is the Hilbert space of all holomorphic functions in $ L^{2}(\mathbb{B}_{n},\mathrm{d}v_{\alpha}) $, where
$$ \mathrm{d}v_{\alpha}(w) = c_{\alpha}(1-|w|^{2})^{\alpha} \mathrm{d}v(w) \ , $$
with $c_{\alpha} = \frac{\Gamma(\alpha+1+n)}{n!\Gamma(\alpha+1)}$ a normalizing constant and $\mathrm{d}v$ standing for the normalized Lebesgue measure on $\mathbb{B}_{n}$. Moreover, $A^{2}_{\alpha}(\mathbb{B}_{n})$ is a reproducing kernel Hilbert space with the inner product inherited from $ L^{2}(\mathbb{B}_{n},\mathrm{d}v_{\alpha}) $ and with kernel $$ K^{\alpha}(z,w) = \frac{1}{(1-\langle z,w \rangle)^{n+1+\alpha}} \ , \ z,w \in \mathbb{B}_{n} \ .$$
When $\alpha = 0$, $\mathrm{d}v_{0} = \mathrm{d}v $ and we write $ A^{2}(\mathbb{B}_{n}) $ and $K(z,w)$ instead of $ A^{2}_{0}(\mathbb{B}_{n}) $ and $K^{0}(z,w)$, respectively. 
The function $K^{\alpha}$ is known as the (weighted) Bergman kernel and it satisfies the reproducing property
\begin{eqnarray}
f(z) 
& = & \int_{\mathbb{B}_{n}} K^{\alpha}(z,w)f(w) \ \mathrm{d}v_{\alpha}(w) \nonumber \ ,
\end{eqnarray}
for any $f \in A^{2}_{\alpha}(\mathbb{B}_{n})$, $z \in \mathbb{B}_{n}$. 
This integral expression is used to define Toeplitz operators. For a bounded Lebesgue measurable funtion $a$ or a regular Borel measure $\mu$, both on $\mathbb{B}_{n}$, their Toeplitz operators $\mathcal{T}_{a}$ and $\mathcal{T}_{\mu}$ are defined by
\begin{equation}\label{definicionoperadorclasico}
\mathcal{T}_{a}f(z) = \int_{\mathbb{B}_{n}} K^{\alpha}(z,w)a(w)f(w) \ \mathrm{d}v_{\alpha}(w) \ ,
\end{equation}
and
\begin{equation}\label{definicionoperadorconmedida}
\mathcal{T}_{\mu}f(z) = c_{\alpha}\int_{\mathbb{B}_{n}} K^{\alpha}(z,w)f(w) (1-|w|^{2})^{\alpha} \ \mathrm{d}\mu(w) \ ,
\end{equation}
for any $f \in  A^{2}(\mathbb{B}_{n}) $, $ z \in \mathbb{B}_{n} $. In each case, $a$ or $\mu$ is called the symbol of the respective Toeplitz operator. Definition (\ref{definicionoperadorclasico}) is the classical for a Toeplitz operator and is equivalent to the compression to $A^{2}_{\alpha}(\mathbb{B}_{n})$ of the multiplying operator $M_{a}f=af$, while $\mathcal{T}_{\mu}$ is a natural generalization of $\mathcal{T}_{a}$, as $\mathcal{T}_{a}$ is the special case when $\mu$ is absolutely continuous with respect to the Lebesgue measure and $\mathrm{d}\mu = a \mathrm{d}v$.

We present the class of Toeplitz operators we will work with.
\begin{definition}\label{definiciontoeplitz}
Let $\Gamma \subset \mathbb{B}_{n}$ be 
a Riemannian submanifold of dimension $d$, $ \mathrm{d}\sigma $ the Riemannian volume element on $\Gamma$, and $ a: \Gamma \rightarrow \mathbb{C} $ a function such that $a \in C^{\infty}_{0}(\Gamma)$. We define the Toeplitz operator $T_{a \mathrm{d}\sigma}$ by
\begin{equation}\label{definiciontoeplitzeq}
T_{a \mathrm{d}\sigma}f(z) = c_{\alpha} \int_{\Gamma} K^{\alpha}(z,w)f(w)a(w)(1-|w|^{2})^{\alpha} \ \mathrm{d}\sigma(w) \ ,
\end{equation}
for every $ f \in A_{\alpha}^{2}(\mathbb{B}_{n}) $, $z \in \mathbb{B}_{n}$.
\end{definition}

\begin{myremark}
We would like to recall that $\Gamma \subset \mathbb{B}_{n}$ being a Riemannian submanifold means that the inclusion 
$\iota: \Gamma \rightarrow \mathbb{B}_{n}$ is a smooth embedding (an immersion which is also a topological embedding) and that the Riemannian metric on $\Gamma$ is the pullback of the Riemannian metric on $\mathbb{B}_{n}$ (See for example \cite{lee}, p.15,132). The Riemannian metric on $\mathbb{B}_{n}$ will always be taken to be the one obtained from the Bergman Hermitian metric (see further below), which makes $\mathbb{B}_{n}$ a model for hyperbolic geometry.
\end{myremark}

\begin{myremark}
\begin{itemize}
\item[•]
Although not made explicit in the notation, $T_{a\mathrm{d}\sigma}$ depends on $\alpha$. We are interested in the asymptotic behavior of $T_{a\mathrm{d}\sigma}$ as $\alpha$ tends to infinity.
\item[•]
If $d=2n$, $\Gamma$ is an open subset of $\mathbb{B}_{n}$ and the measure $ \mathrm{d}\sigma(w) = \frac{\pi^{n}}{n!} \frac{\mathrm{d}v(w)}{(1-|w|^{2})^{n+1}}$ is a multiple of the invariant measure on $\mathbb{B}_{n}$. Therefore $T_{a} := T_{a\mathrm{d}\sigma}$ is the {\it ``classical"} Toeplitz operator (\ref{definicionoperadorclasico}) with symbol
$ \frac{\pi^{n}}{n!} \frac{a(w)}{(1-|w|^{2})^{n+1}} \ . $ Explicitly $$ \mathcal{T}_{a} = T_{\scriptscriptstyle \frac{n!}{\pi^{n}} a(\cdot)(1-|\cdot|^{2})^{n+1}} \ . $$
\item[•]
If $d<2n$, the resulting symbol measure $\mu = a \mathrm{d}\sigma$ is singular with respect to the Lebesgue measure on $\mathbb{B}_{n}$. Therefore, $T_{a\mathrm{d}\sigma}$ is of the form (\ref{definicionoperadorconmedida}) and can not be written as (\ref{definicionoperadorclasico}).
\item[•] As $a$ has compact support, $T_{a\mathrm{d}\sigma}$ is bounded, compact, and belongs to each Schatten class $S_{p}(A_{\alpha}^{2}(\mathbb{B}_{n}))$, $p>0$.
\item[•] Condition $a \in C^{\infty}(\Gamma)$ can be changed to $a \in L^{\infty}(\Gamma)$ without affecting the definition of $T_{a\mathrm{d}\sigma}$. In fact, it is not needed for a substantial part of the results. However, it is included in the definition since it is necessary for the main theorem.
\end{itemize}
\end{myremark}

Prior to stating the Szeg\"{o} limit theorem we need to introduce some geometric definitions (see, for example, \cite{huybrechts}, \cite{krantz}, \cite{libermann}).

By writing the local holomorphic coordinates as $ z_{\ell}=x_{\ell}+i y_{\ell} $, $\ell=1, \sdots , n$, $\mathbb{B}_{n}$ can be naturally regarded both as a complex manifold of dimension $n$ and as a smooth manifold of dimension $2n$. For $p \in \mathbb{B}_{n}$ the respective real and complex tangent spaces are given in terms of basic tangent vectors by
$$ \mathrm{T}_{p}\mathbb{B}_{n} = \mathrm{span}\left\{ \tngvec[x]{\ell}{p}, \tngvec[y]{\ell}{p} \right\}_{\scriptscriptstyle \ell=1, \dots , n} \mbox{and} \ \ \mathrm{T}_{p}\mathbb{B}_{n}^{\mathbb{C}} = \mathrm{span}\left\{ \tngvec[z]{\ell}{p} ,\tngvec[\overline{z}]{\ell}{p} \right\}_{\scriptscriptstyle \ell=1, \dots , n} \ , $$
with the relation between the elements of these bases given by the Wirtinger derivatives 
\begin{equation}\label{baseconwirtingerderivatives}
 \frac{\partial}{\partial z_{\ell}}  = \frac{1}{2} \left( \frac{\partial}{\partial x_{\ell}} - i \frac{\partial}{\partial y_{\ell}} \right) \  \ , \ \ \frac{\partial}{\partial \overline{z}_{\ell}} = \frac{1}{2} \left( \frac{\partial}{\partial x_{\ell}} + i \frac{\partial}{\partial y_{\ell}} \right) \ . \nonumber
\end{equation}

The Hermitian metric $b$ defined by $$ b_{jk}(p) := \frac{1}{n+1} \frac{\partial^{2}}{\partial z_{j} \partial \overline{z}_{k}} \log K(z,z) \Big\rvert_{z=p}  \ , $$ $ j,k = 1, \sdots , n$ on $\mathrm{T}_{p}\mathbb{B}_{n}^{\mathbb{C}} $ for each $p \in \mathbb{B}_{n}$ is called the Bergman Hermitian metric and it makes $\mathbb{B}_{n}$ a Hermitian manifold. It induces a Riemannian metric $ \mathring{b} $ on the underlying smooth manifold, thus giving it the structure of a Riemannian manifold. The metric $ \mathring{b} $ is described in terms of $b$ by the matrix
\begin{equation}
\mathring{B}(p) := \left( \mathring{b}_{jk}(p) \right) =
\begin{pmatrix}
\phantom{-} \mathrm{Re} \ b_{jk}(p) & \phantom{-} \mathrm{Im} \ b_{jk}(p) \\
-\mathrm{Im} \ b_{jk}(p) & \phantom{-} \mathrm{Re} \ b_{jk}(p)
\end{pmatrix}  \ . \nonumber
\end{equation}
The endomorphism $J$ defined by 
$$ J\left( \tngvec[x]{\ell}{p} \right) = \tngvec[y]{\ell}{p} \ , \ \quad J\left( \tngvec[y]{\ell}{p} \right) = -\tngvec[x]{\ell}{p} \  $$
is a complex structure on $\mathrm{T}_{p}\mathbb{B}_{n}$ that is compatible with $\mathring{b}$ in the sense that $\mathring{b}$ is invariant under $J$. This complex structure is in turn induced by the natural one of $\mathbb{R}^{2n}$,
$ (x_{1}, \sdots x_{n}, y_{1}, \sdots , y_{n}) \mapsto (-y_{1}, \sdots , -y_{n}, x_{1}, \sdots , x_{n})$. We also consider the symplectic forms defined by
$$ \omega_{p}(\xi , \nu):= \mathring{b}_{p} \big( \xi , J \nu \big) \  $$
for $\xi , \nu \in \mathrm{T}_{p}\mathbb{B}_{n}$. If $V \subseteq \mathrm{T}_{p}\mathbb{B}_{n}$ is a subspace, we denote by
\begin{equation}\label{definiciondelcomplementosimplectico}
    V^{\omega} := \big\{ \xi \in \mathrm{T}_{p}\mathbb{B}_{n} \ : \ \omega_{p}(\xi , v) = 0 \ \ \forall v \in V \big\} 
\end{equation}
its symplectic complement.
\begin{definition}\label{definiciondeisotropicoyco-isotropico}
We say that $V$ is isotropic (respectively co-isotropic, Lagrangian) if $ V \subseteq V^{\omega} $ (respectively $V \supseteq V^{\omega} $, $ V = V^{\omega} $).
\end{definition}

\begin{example}\label{ejemplodeespaciosisosycoisos}
From the skew-symmetry of $\omega_{p}$ it follows that any subspace of dimension one is isotropic and any of co-dimension one is co-isotropic. In addition, the total subspace $\mathrm{T}_{p}\mathbb{B}_{n}$ is trivially co-isotropic.
\end{example}

For $\Gamma \subseteq \mathbb{B}_{n}$ a submanifold, let $\mathrm{D}_{p}\iota: \mathrm{T}_{p}\Gamma \rightarrow \mathrm{T}_{p}\mathbb{B}_{n} $ be the differential of the inclusion map $\iota:\Gamma \rightarrow \mathbb{B}_{n} $ at the point $p\in \Gamma$. Then $\mathrm{D}_{p}\iota \big( \mathrm{T}_{p}\Gamma \big)$ is a subspace of $\mathrm{T}_{p}\mathbb{B}_{n}$ of dimension $d$.

\begin{definition}
Let $\Gamma \subseteq \mathbb{B}_{n}$ be a submanifold. We say that $\Gamma$ is isotropic (respectively co-isotropic, Lagrangian) if for every $p \in \Gamma$, the subspace $\mathrm{D}_{p}\iota \big( \mathrm{T}_{p}\Gamma \big) \subseteq \mathrm{T}_{p}\mathbb{B}_{n} $ is isotropic (respectively co-isotropic, Lagrangian).  
\end{definition}

\begin{example}
Following Example \ref{ejemplodeespaciosisosycoisos}, any curve is isotropic and any hypersurface is co-isotropic. Also, $\mathbb{B}_{n}$ or any open subset of $\mathbb{B}_{n}$ is co-isotropic.
\end{example}

\subsection{Statement of the main result}

Our Szeg\"{o} limit theorem will apply to submanifolds that are isotropic or co-isotropic. In order to include both cases in the same statement, we introduce the notation
\begin{equation}
d' := \left\{ \begin{array}{cl}
d & \Gamma \mbox{ \ isotropic} \\
2n-d & \Gamma \mbox{ \ co-isotropic}
\end{array} \right. \ . \nonumber
\end{equation}

We also normalize $T_{a \mathrm{d}\sigma}$ by means of
\begin{equation}\label{lanormalizacionprimeravez}
\widehat{T}_{a\mathrm{d}\sigma} := \frac{n!}{2^{\frac{d'}{2}}\pi^{\frac{d}{2}}}\alpha^{-n+\frac{d}{2}} T_{a\mathrm{d}\sigma} \ .
\end{equation}
We will prove in Corollary \ref{corolariodelanormalizacion} that $ \big\vert \big\vert \widehat{T}_{a\mathrm{d}\sigma} \big\vert \big\vert_{\mathcal{L}\left( A^{2}_{\alpha}(\mathbb{B}_{n}) \right)} = O(1) $ as $\alpha$ tends to infinity, thus justifying the usage of the term normalization.

We can now state the main result of this paper. It is a Bergman space analogue of the main result of \cite{SalvadorAlejandroUribe}.

\begin{theorem}\label{teoremaszego}{\bf [Szeg\"{o} limit theorem]}
Let $\Gamma \subseteq \mathbb{B}_{n}$ be a isotropic or co-isotropic submanifold, $a \in C^{\infty}_{0}(\Gamma)$ a positive function, and $R>0$ such that $\sup\limits_{\xi \in \Gamma} \frac{a(\xi)}{(1-|\xi|^{2})^{n+1}}< R$. If $\varphi$ is a function such that $\frac{\varphi(t)}{t^{p}} \in C[0,R]$ for some $p>0$, then
\begin{equation}\label{eqteoremaszego}
\lim\limits_{\alpha \rightarrow \infty} \left( \frac{\pi}{\alpha} \right)^{\frac{d}{2}} \mathrm{Tr}\left( \varphi\left(\widehat{T}_{a\mathrm{d}\sigma}\right) \right) = \frac{1}{2^{\frac{d'}{2}}} \int_{\Gamma} \quessianobb{\varphi} \ \mathrm{d}\sigma(\xi) \ ,
\end{equation}
where $\mathcal{Q}_{\varepsilon}$ is the operator
$$ \mathcal{Q}_{\varepsilon}(\varphi)(t) := \frac{1}{\Gamma(\varepsilon)} \int_{0}^{t} \varphi(s) \left( \ln \frac{t}{s} \right)^{\varepsilon-1} \ \frac{\mathrm{d}s}{s} $$
when $\varepsilon > 0$ ($\Gamma(\varepsilon)$ is the Gamma function evaluated at $\varepsilon$) and the identity when $\varepsilon = 0$.
\end{theorem}


\begin{example}[The absolutely continuous case]
The case $\Gamma=\mathbb{B}_{n}$ is trivially co-isotropic
, with (\ref{eqteoremaszego}) becoming
\begin{equation}
\lim\limits_{\alpha \rightarrow \infty} \left( \frac{\pi}{\alpha} \right)^{n} \mathrm{Tr}\left( \varphi\left(\widehat{T}_{a}\right) \right) = \frac{\pi^{n}}{n!} \int_{\mathbb{B}_{n}} \varphi\left( \frac{a(w)}{(1-|w|^{2})^{n+1}} \right) \ \frac{\mathrm{d}v(w)}{(1-|w|^{2})^{n+1}} \ . \nonumber
\end{equation}
In terms of the {\it ``classical"} operator $\mathcal{T}_{a}$ (\ref{definicionoperadorclasico}), we have $ \mathcal{T}_{a} = \widehat{T}_{\scriptscriptstyle a(\cdot)(1-|\cdot|^{2})^{n+1}} $ and so
\begin{equation}
\lim\limits_{\alpha \rightarrow \infty} \frac{n!}{\alpha^{n}} \mathrm{Tr}\big( \varphi\left( \mathcal{T}_{a} \right) \big) = \int_{\mathbb{B}_{n}} \varphi\left( a(w) \right) \ \frac{\mathrm{d}v(w)}{(1-|w|^{2})^{n+1}} \ . \nonumber
\end{equation}
The case $n=1$ appeared recently as Theorem 1.1 in \cite{elnuevo} without the smoothness hypothesis for the symbol $a$, but for a narrower family of functions $\varphi$.
\end{example}

The article has the following structure. In the following section, we describe some aspects of the intrinsic geometry of $\Gamma$. We obtain a relation involving the Riemannian metric and the symplectic form in $\Gamma$ that will be
needed later. In Section \ref{laseccionconlatrazadelacomposicion} we derive an integral representation for the trace of a composition of multiple Toeplitz operators with symbols sharing the same support and then apply the method of stationary phase to describe the asymptotic behavior of this integral as $\alpha$ grows to infinity. The expression we obtain relies on computing the determinant of a given Hessian matrix, which we leave for Appendix \ref{laseccionconloscalculosparaeldeterminante}. The proof of the Szeg\"{o} limit theorem is delivered in Section \ref{laseccionconelteorema}, first for polynomials and then for a wider family of continuous functions. A key step in the proof involves bounding the traces for small powers (exponents $0<p<1$) of the Toeplitz operators. This is done in Appendix \ref{seccioncondemodelaclaimdepchica}. Section \ref{laseccionconloscorolarios} is reserved for presenting some corollaries to the main theorem regarding eigenvalue distribution and Schatten norm estimates. We also expose an example illustrating the given results.

The line of reasoning used here follows that in \cite{SalvadorAlejandroUribe}, with some of the proofs presented being adaptations to the Bergman setting of the ones therein. In particular, some proofs in Section \ref{laseccionconloscorolarios} and Appendix \ref{subseccionochopuntouno} proceed virtually unmodified and are therefore omitted for the sake of briefness.

Throughout the article, we will use the lone capital letter $C$ to represent an arbitrary constant that may change from line to line. When deemed appropriate, we will highlight its independence from other variables. The symbol $B(0,\varepsilon)$ will denote the euclidean ball in $\mathbb{R}^{d}$ with radius $\varepsilon>0$ centered at the origin.

\section{The Riemannian metric on \texorpdfstring{$\Gamma$}{Γ}}\label{lasecciondelametricariemannianadeGamma}

Recall that for each $ v \in \mathrm{T}_{p}\Gamma $, 
the differential of the inclusion $\iota: \Gamma \rightarrow \mathbb{B}_{n}$ induces a tangent vector $ \mathrm{D}_{p}\iota (v)$ to $\mathbb{B}_{n}$ at $p=\iota(p)$. As a submanifold $\Gamma \subseteq \mathbb{B}_{n}$, the Riemannian metric $g$ in $\Gamma$ is given by the pullback of the metric $\mathring{b}$ in $\mathbb{B}_{n}$. More precisely
\[ g_{p}( u,v )  = \mathring{b}_{p}\big( \mathrm{D}_{p}\iota (u) , \mathrm{D}_{p}\iota (v) \big) \ , \] 
$\forall \ u , v \ \in \mathrm{T}_{p}\Gamma$. We also define a skew-symmetric, nondegenerate bilinear form $h$ in terms of the symplectic form $\omega_{p}$ by
\[ h_{p}( u,v )  = \omega_{p}\big( \mathrm{D}_{p}\iota (u) , \mathrm{D}_{p}\iota (v) \big) \ , \] 
$\forall \ u , v \ \in \mathrm{T}_{p}\Gamma$.\par

Let $\gamma : B(0,1) \rightarrow \Gamma$ be a parametrization, so that for each ${p=\gamma(t)\in \Gamma}$, 
the tangent space $\mathrm{T}_{p}\Gamma$ is spanned by the tangent vectors $\tngvec{j}{t} $, $j=1,\sdots , d$. Denote $\gamma_{\ell}:B(0,1) \rightarrow \mathbb{C}$, $\ell = 1, \dots , n$ its component functions. The following lemma relating $ g_{p} $ and $h_{p}$ looks rather technical, but will be key in Section \ref{laseccionconelcomportamientoasintotico}.

\begin{lemma}\label{lemaconlasrealcionesdelasformasbilinealesengamma}
When taking into account the basis $ \left\{ \tngvec{j}{t} \right\}_{j=1,\dots ,d} $ of $\mathrm{T}_{p}\Gamma$, the components of the forms $g_{p}$ and $h_{p}$ satisfy the relations
\begin{eqnarray}
g_{jk}(p)+ih_{jk}(p) & = & \sum_{\ell , r=1}^{n} b_{\ell r}(p) \cdot \partial_{j}\gamma_{\ell}(t) \cdot \partial_{k}\overline{\gamma_{r}}(t) \ , \label{relaciondelasmetricasengamma} \\
& \ & \ \nonumber \\
g_{jk}(p)-ih_{jk}(p) & = & \sum_{\ell , r=1}^{n} \overline{b_{\ell r}(p)} \cdot \partial_{j}\overline{\gamma_{\ell}}(t) \cdot \partial_{k}\gamma_{r}(t) \ , \label{relaciondelasmetricasengammaconjugada}
\end{eqnarray}
for $j,k=1,\sdots ,d$.
\end{lemma}
\begin{proof}
For convenience, we will write $ \tngvec[x]{n+\ell}{p} = \tngvec[y]{\ell}{p} $. First we note that for each $j=1,\sdots , d$
$$ \mathrm{D}_{p}\iota \left( \tngvec{j}{t} \right) = \sum_{\ell=1}^{2n} \partial_{j}(\iota \circ \gamma)_{\ell}(t) \tngvec[x]{\ell}{p} \ , $$
where $ (\iota \circ \gamma)_{\ell} $ denotes the $\ell$-th component function of $\iota \circ \gamma :B(0,1) \subset \mathbb{R}^{d} \rightarrow \mathbb{B}_{n} \subseteq \mathbb{R}^{2n} $, and $\partial_{j}(\iota \circ \gamma)_{\ell}(t)$ its partial derivative respect to the $j$-th variable valued at $t$. This is obtained straight from the definitions and using the chain rule. We may now compute

{\allowdisplaybreaks
\begin{align}
    g_{jk}(p) & = \mathring{b}_{p}\left( \mathrm{D}_{p}\iota \left( \tngvec{j}{t} \right) , \mathrm{D}_{p}\iota \left( \tngvec{k}{t} \right) \right) \nonumber \\
    & = \sum_{\ell , r = 1}^{2n} \partial_{j}(\iota \circ \gamma)_{\ell}(t) \partial_{k}(\iota \circ \gamma)_{r}(t) \mathring{b}_{p}\left( \tngvec[x]{\ell}{p} , \tngvec[x]{r}{p} \right) \nonumber \\
    & = \sum_{\ell , r = 1}^{2n} \partial_{j}(\iota \circ \gamma)_{\ell}(t) \partial_{k}(\iota \circ \gamma)_{r}(t) \mathring{b}_{\ell r}(\gamma(t)) \nonumber \\  
    & = \sum_{\ell , r = 1}^{n} \bigg[ \bigg( \partial_{j}(\iota \circ \gamma)_{\ell}(t) \partial_{k}(\iota \circ \gamma)_{r}(t) + \partial_{j}(\iota \circ \gamma)_{n+\ell}(t) \partial_{k}(\iota \circ \gamma)_{n+r}(t) \bigg) \nonumber \\
    & \quad \qquad \cdot \mathrm{Re } b_{\ell r}(\gamma(t)) + \bigg( \partial_{j}(\iota \circ \gamma)_{\ell}(t) \partial_{k}(\iota \circ \gamma)_{n+r}(t) \nonumber \\
    & \quad \qquad \qquad - \partial_{j}(\iota \circ \gamma)_{n+\ell}(t) \partial_{k}(\iota \circ \gamma)_{r}(t) \bigg) \mathrm{Im } b_{\ell r}(\gamma(t)) \bigg] \ . \nonumber  
\end{align}
}

In a similar fashion, we get

\begin{align}
    h_{jk}(p) & = \mathring{b}_{p}\left( \mathrm{D}_{p}\iota \left( \tngvec{j}{t} \right) , J \left( \mathrm{D}_{p}\iota \left( \tngvec{k}{t} \right) \right) \right) \nonumber \\
    & = \sum_{\ell , r = 1}^{n} \bigg[ \bigg( \partial_{j}(\iota \circ \gamma)_{\ell}(t) \partial_{k}(\iota \circ \gamma)_{r}(t) + \partial_{j}(\iota \circ \gamma)_{n+\ell}(t) \partial_{k}(\iota \circ \gamma)_{n+r}(t) \bigg) \nonumber \\
    & \quad \qquad \cdot \mathrm{Im } b_{\ell r}(\gamma(t)) + \bigg( \partial_{j}(\iota \circ \gamma)_{n+\ell}(t) \partial_{k}(\iota \circ \gamma)_{r}(t) \nonumber \\
    & \quad \qquad \qquad - \partial_{j}(\iota \circ \gamma)_{\ell}(t) \partial_{k}(\iota \circ \gamma)_{n+r}(t) \bigg) \mathrm{Re } b_{\ell r}(\gamma(t)) \bigg] \ . \nonumber
\end{align}

Using these two relations

\begin{align}
   g_{jk}(p)+ih_{jk}(p) & {\scriptstyle =} \sum_{\ell , r=1}^{n} \bigg[ \bigg( \partial_{j}(\iota \circ \gamma)_{\ell} \partial_{k}(\iota \circ \gamma)_{r} + \partial_{j}(\iota \circ \gamma)_{n+\ell} \partial_{k}(\iota \circ \gamma)_{n+r} \bigg) \nonumber \\
   & \quad \qquad \cdot \bigg( \mathrm{Re } b_{\ell r}(p) + i \mathrm{Im } b_{\ell r}(p) \bigg) + \bigg( \partial_{j}(\iota \circ \gamma)_{n+\ell} \partial_{k}(\iota \circ \gamma)_{r} \nonumber \\
   & \quad \qquad \quad - \partial_{j}(\iota \circ \gamma)_{\ell} \partial_{k}(\iota \circ \gamma)_{n+r} \bigg) \bigg( i \mathrm{Re } b_{\ell r}(p) - \mathrm{Im } b_{\ell r}(p) \bigg) \bigg] \nonumber \\
   & {\scriptstyle =} \sum_{\ell , r=1}^{n} b_{\ell r}(p) \bigg[ \bigg( \partial_{j}(\iota \circ \gamma)_{\scriptstyle \ell} \partial_{k}(\iota \circ \gamma)_{\scriptstyle r} + \partial_{j}(\iota \circ \gamma)_{ n+\ell} \partial_{k}(\iota \circ \gamma)_{ n+r} \bigg) \nonumber \\
   & \quad \qquad + i \bigg( \partial_{j}(\iota \circ \gamma)_{n+\ell} \partial_{k}(\iota \circ \gamma)_{r} - \partial_{j}(\iota \circ \gamma)_{\ell} \partial_{k}(\iota \circ \gamma)_{n+r} \bigg) \bigg] \nonumber \\
   & {\scriptstyle =} \sum_{\ell , r=1}^{n} b_{\ell r}(p) \bigg[ \partial_{j} \Big( ( \iota \circ \gamma)_{\ell} + i (\iota \circ \gamma)_{n+\ell} \Big) \nonumber \\
   & \quad \qquad \cdot \partial_{k} \Big( (\iota \circ \gamma)_{r} - i (\iota \circ \gamma)_{n+r} \Big) \bigg] \nonumber \\
   & {\scriptstyle =} \sum_{\ell , r=1}^{n} b_{\ell r}(p) \cdot \partial_{j}\gamma_{\ell}(t) \cdot \partial_{k}\overline{\gamma_{r}}(t) \ , \nonumber
\end{align}

from where (\ref{relaciondelasmetricasengammaconjugada}) is also deduced.

Notice that, as we have done when writing the local holomorphic coordinates of $ \mathbb{B}_{n}$, we have used $ (\iota \circ \gamma)_{\ell} + {i (\iota \circ \gamma)_{n+\ell}} = \gamma_{\ell} $.
\end{proof}

\section{The trace of a multiple composition}
\label{laseccionconlatrazadelacomposicion}

The first step towards a proof of Theorem \ref{teoremaszego} relies on representing in terms of an integral the composition of multiple Toeplitz operators with symbols sharing the same support.

\subsection{An integral expression for the composition}

We recall some notions from the classic theory of Berezin symbols and kernels (\cite{berezin}). Let $A \in \mathcal{L}\big(A^{2}_{\alpha}(\mathbb{B}_{n})\big)$. The associated kernel $K^{\alpha}_{A}$ is the function
$$ K^{\alpha}_{A}(z,w) := Ak_{w}(z) = \langle Ak_{w} , k_{z} \rangle_{A^{2}_{\alpha}(\mathbb{B}_{n})} \ ,$$
where $k_{z}= \frac{K^{\alpha}(\cdot,z)}{\sqrt{K^{\alpha}(z,z)}} = \left( \frac{\sqrt{1-|z|^{2}}}{1-\langle \cdot,z \rangle} \right)^{n+1+\alpha}$ denotes the normalized reproducing kernel at $z$.
In the case when $A$ belongs to the trace class, then $K^{\alpha}_{A}$ satisfies
\begin{equation}\label{nucleoasociado.traza}
\mathrm{Tr}(A) = \int_{\mathbb{B}_{n}} K_{A}^{\alpha}(z,z) \ \mathrm{d}v_{\alpha}(z) \ ,
\end{equation}
and whenever $B$ is another bounded operator we have
\begin{equation}\label{nucleoasociado.dos}
K^{\alpha}_{AB}(z,w) = \int_{\mathbb{B}_{n}} K_{A}^{\alpha}(z,y) K_{B}^{\alpha}(y,w) \ \mathrm{d}v_{\alpha}(y) \ .
\end{equation}

For the instance when $A=T_{a\mathrm{d}\sigma}$ is a Toeplitz operator, Fubini's theorem and the reproducing property of the Bergman kernels yield in the expression
\begin{equation}\label{nucleoasociado.toeplitz}
K^{\alpha}_{T_{a\mathrm{d}\sigma}}(z,w) = c_{\alpha} \int_{\Gamma} K^{\alpha}(z,\xi) K^{\alpha}(\xi,w) \ a(\xi)(1-|\xi|^{2})^{\alpha} \ \mathrm{d}\sigma(\xi) \ .
\end{equation}

We present a formula for the trace of a composition of $ m $ Toeplitz operators 

\begin{lemma}\label{1.5formuladelatraza}
Let $m\geq 2$ be a natural number, $ a_{1}, \sdots , a_{m} \in L^{\infty}(\Gamma) $, and ${ S:= T_{a_{1}\mathrm{d}\sigma} \circ \scdots \circ T_{a_{m}\mathrm{d}\sigma} }$. Then

\begin{multline}\label{latrazadelacomposicion}
\mathrm{Tr}(S) = c_{\alpha}^{m} \int\limits_{\Gamma^{m}} \prod_{\scriptscriptstyle j (\mathrm{mod} \ m)} \bigg[ \left( \frac{1-|\xi_{j}|^{2}}{1-\langle \xi_{j} , \xi_{j+1} \rangle} \right)^{\alpha} \frac{a_{j}(\xi_{j})}{(1-\langle \xi_{j} , \xi_{j+1} \rangle )^{n+1}} \bigg] \\
\cdot \mathrm{d}\sigma(\xi_{1}) \cdots \mathrm{d}\sigma(\xi_{m}) \ .
\end{multline}
\end{lemma}
\begin{proof}
By induction on the relation (\ref{nucleoasociado.dos}) one can prove that

\begin{multline}
    K_{S}^{\alpha}(z,w) = \int\limits_{\mathbb{B}_{n}^{m-1}} K^{\alpha}_{T_{a_{1}\mathrm{d}\sigma}}(z,y_{1}) \left[ \prod_{j=2}^{m-1} K^{\alpha}_{T_{a_{j}\mathrm{d}\sigma}}(y_{j-1},y_{j}) \right] K^{\alpha}_{T_{a_{m}\mathrm{d}\sigma}}(y_{m-1},w) \\
    \cdot \mathrm{d}v_{\alpha}(y_{m-1}) \cdots \mathrm{d}v_{\alpha}(y_{1}) \ . \nonumber
\end{multline}

We write $ y_{0}=z $ and $y_{m}=w$. Replacing all $ K^{\alpha}_{T_{a_{j}\mathrm{d}\sigma}} $ for their equivalent form in (\ref{nucleoasociado.toeplitz}), and using Fubini's theorem together with the reproducing property of the Bergman kernel, we get

\begin{align}
    K_{S}^{\alpha}(z,w) & = c_{\alpha}^{m} \int\limits_{\mathbb{B}_{n}^{m-1}} \prod_{j=1}^{m} \left[ \int\limits_{\Gamma}K^{\alpha}(y_{j-1},\xi_{j}) K^{\alpha}(\xi_{j},y_{j}) \ a_{j}(\xi_{j})(1-|\xi_{j}|^{2})^{\alpha} \ \mathrm{d}\sigma(\xi_{j}) \right] \nonumber \\
    & \quad \qquad \cdot \mathrm{d}v_{\alpha}(y_{m-1}) \cdots \mathrm{d}v_{\alpha}(y_{1}) \nonumber \\
    & = c_{\alpha}^{m} \int\limits_{\Gamma^{m}} \int\limits_{\mathbb{B}_{n}^{m-1}} \prod_{j=1}^{m} \Big[ K^{\alpha}(y_{j-1},\xi_{j}) K^{\alpha}(\xi_{j},y_{j}) \Big] \mathrm{d}v_{\alpha}(y_{1}) \cdots \mathrm{d}v_{\alpha}(y_{m-1}) \nonumber \\
    & \quad \qquad \cdot \prod_{\ell=1}^{m} \left[ {a_{\ell}(\xi_{\ell})}(1-|\xi_{\ell}|^{2})^{\alpha} \right] \mathrm{d}\sigma(\xi_{1}) \cdots \mathrm{d}\sigma(\xi_{m}) \nonumber \\
    & = c_{\alpha}^{m} \int\limits_{\Gamma^{m}} K^{\alpha}(y_{0},\xi_{1})K^{\alpha}(\xi_{m},y_{m}) \prod_{j=1}^{m-1} \left[ \int\limits_{\mathbb{B}_{n}} K^{\alpha}(\xi_{j},y_{j})k_{\xi_{j+1}}(y_{j}) \ \mathrm{d}v_{\alpha}(y_{j}) \right] \nonumber \\
    & \quad \qquad \cdot \prod_{\ell=1}^{m} \left[ a_{\ell}(\xi_{\ell})(1-|\xi_{\ell}|^{2})^{\alpha} \right] \mathrm{d}\sigma(\xi_{1}) \cdots \mathrm{d}\sigma(\xi_{m}) \nonumber \\
    & = c_{\alpha}^{m} \int\limits_{\Gamma^{m}} K^{\alpha}(z,\xi_{1})K^{\alpha}(\xi_{m},w) \prod_{j=1}^{m-1} \Big[ K^{\alpha}(\xi_{j},\xi_{j+1}) \Big] \nonumber \\
    & \quad \qquad \cdot \prod_{\ell=1}^{m} \left[ a_{\ell}(\xi_{\ell})(1-|\xi_{\ell}|^{2})^{\alpha} \right] \mathrm{d}\sigma(\xi_{1}) \cdots \mathrm{d}\sigma(\xi_{m})  \ . \nonumber
\end{align}

Using the trace formula in (\ref{nucleoasociado.traza}),

\begin{align}
    \mathrm{Tr}(S) & = \int_{\mathbb{B}_{n}} K_{S}^{\alpha}(z,z) \ \mathrm{d}v_{\alpha}(z) \nonumber \\
    & = c_{\alpha}^{m} \int\limits_{\mathbb{B}_{n}} \int\limits_{\Gamma^{m}} K^{\alpha}(z,\xi_{1})K^{\alpha}(\xi_{m},z) \prod_{j=1}^{m-1} \Big[ K^{\alpha}(\xi_{j},\xi_{j+1}) \Big] \nonumber \\
    & \quad \qquad \cdot \prod_{\ell=1}^{m} \left[ a_{\ell}(\xi_{\ell})(1-|\xi_{\ell}|^{2})^{\alpha} \right] \mathrm{d}\sigma(\xi_{1}) \cdots \mathrm{d}\sigma(\xi_{m}) \ \mathrm{d}v_{\alpha}(z) \nonumber \\
    & = c_{\alpha}^{m} \int\limits_{\Gamma^{m}} \int\limits_{\mathbb{B}_{n}} K^{\alpha}(\xi_{m},z)k_{\xi_{1}}(z) \ \mathrm{d}v_{\alpha}(z) \prod_{j=1}^{m-1} \Big[ K^{\alpha}(\xi_{j},\xi_{j+1}) \Big] \nonumber \\
    & \quad \qquad \cdot \prod_{\ell=1}^{m} \left[ a_{\ell}(\xi_{\ell})(1-|\xi_{\ell}|^{2})^{\alpha} \right] \mathrm{d}\sigma(\xi_{1}) \cdots \mathrm{d}\sigma(\xi_{m}) \nonumber \\
    & = c_{\alpha}^{m} \int\limits_{\Gamma^{m}} \prod_{j (\mathrm{mod} m)} \Big[ K^{\alpha}(\xi_{j},\xi_{j+1}) \ a_{j}(\xi_{j})(1-|\xi_{j}|^{2})^{\alpha} \Big] \ \mathrm{d}\sigma(\xi_{1}) \cdots \mathrm{d}\sigma(\xi_{m}) \ , \nonumber
\end{align}

from where the statement follows.
\end{proof}

We rewrite the expression (\ref{latrazadelacomposicion}) in the form
\begin{multline}\label{formulatraza}
\mathrm{Tr}(S) = c_{\alpha}^{m} \int\limits_{\Gamma^{m}} \mathrm{exp}\left[ i\alpha \left( i \Log{\prod_{j=1}^{m}} \frac{1 {\scriptstyle - }\langle \xi_{j} , \xi_{j+1} \rangle}{1 {\scriptstyle - } |\xi_{j}|^{2}} \right) \right] \left( \prod_{j=1}^{m} \frac{a_{j}(\xi_{j})}{(1 - \langle \xi_{j} , \xi_{j+1} \rangle)^{n+1}} \right) \\ 
\cdot \mathrm{d}\sigma(\xi_{1}) \cdots \mathrm{d}\sigma(\xi_{m}) \ ,
\end{multline}
with the understanding that the subindexes are always taken $\mathrm{mod} \ m$ (i.e. $ \xi_{m+1}=\xi_{1} $). We also define
\begin{eqnarray}
\Phi( \xi_{1}, \sdots , \xi_{m} ) & := &  i \Log \prod_{j=1}^{m} \frac{1-\langle \xi_{j} , \xi_{j+1} \rangle}{1-|\xi_{j}|^{2}} \nonumber \\
& = & i \sum_{j=1}^{m} \Log{ \frac{1-\langle \xi_{j} , \xi_{j+1} \rangle}{1-|\xi_{j}|^{2}} }  \label{laformadelaphifaseenterminosdesumadelogs} \\
& = & \frac{i}{n+1} \sum_{j=1}^{m} \big[ \Log K(\xi_{j},\xi_{j}) - \Log K(\xi_{j},\xi_{j+1}) \big] \ \label{laformadelaphifaseenterminosdellogdelak} .
\end{eqnarray}
This notation will be useful when using stationary phase method later.

\begin{remark}\label{remarktraza}
Notice that $\frac{1-\langle \xi_{j} , \xi_{j+1} \rangle}{1-|\xi_{j}|^{2}} \in (-\infty , 0] $ if and only if $\langle \xi_{j} , \xi_{j+1} \rangle \in [1, \infty) $, which does not occur since $ \xi_{j}, \xi_{j+1} \in \mathbb{B}_{n} $.
It follows from (\ref{laformadelaphifaseenterminosdesumadelogs}) that $ \Phi $ is well defined and continuous on $ \Gamma^{m} $.
Moreover, from (\ref{formulatraza}) we find that all the $\xi_{j}$ involved in the integration are the ones in the support of $a_{j}$. If each $a_{j}$ has compact support in $\Gamma$, we conclude that there is an open neighborhood of $(-\infty ,0] $ such that only values outside this neighborhood are used to evaluate $\Log$ in the expression for $\Phi$ when computing (\ref{formulatraza}).
\end{remark}

\subsubsection{The behavior on the diagonal}

We describe the behavior of $ \mathrm{Im}\big( \Phi(\xi_{1}, \sdots , \xi_{m}) \big)$ with the following result. From it we deduce that the critical behavior of $\Phi$ occurs on the diagonal of $\Gamma^{m}$.

\begin{lemma}\label{lemmadeladesigualdaddelafase}
$ \mathrm{Im}\big( \Phi(\xi_{1}, \sdots , \xi_{m}) \big) \geq 0$
for any $(\xi_{1}, \sdots , \xi_{m}) \in \Gamma^{m}$
, with equality if and only if $ \xi_{1} = \xi_{2} = \scdots = \xi_{m} $.
\end{lemma}
We will make use of the following claim.
\begin{description}
\item \begin{claim}\label{proposiciondelasetiquetas}
Let $m \in \mathbb{N}$, $m \geq 2$. The for any $ \{ d_{1}, \sdots , d_{m} \} \subset (0,1)$,
\begin{equation}\label{afirmaciondelasetiquetasprimera}
\prod_{j=1}^{m} \frac{1-d_{j}d_{j+1}}{1-d_{j}^{2}} \geq 1 \ , 
\end{equation}
where $ d_{m+1}:=d_{1} $. Moreover, the equality holds if and only if $ d_{1}=d_{2}= \scdots = d_{m} $.
\end{claim}
\begin{proof}
The proof is straightforward when $d_{1} = d_{2} = \scdots = d_{m}$. Suppose that not all $d_{j}$ are the same and let $d_{0}:=d_{m}$. We must show that the inequality in (\ref{afirmaciondelasetiquetasprimera}) is strict. Notice that if $k'$ is such that $d_{k'+1}=d_{k'}$ then
$$ \frac{1-d_{k'-1}d_{k'}}{1-d_{k'-1}^{2}} \cdot \frac{1-d_{k'}d_{k'+1}}{1-d_{k'}^{2}} = \frac{1-d_{k'-1}d_{k'+1}}{1-d_{k'-1}^{2}} \ , $$
and so 
\begin{eqnarray}
\prod_{j=1}^{m} \frac{1-d_{j}d_{j+1}}{1-d_{j}^{2}} & = & \frac{1-d_{k'-1}d_{k'+1}}{1-d_{k'-1}^{2}} \prod_{\substack{j=1 \\ j \neq k'-1,k'}}^{m} \frac{1-d_{j}d_{j+1}}{1-d_{j}^{2}}  \ . \nonumber
\end{eqnarray}
The rightmost term is a product of the form of (\ref{afirmaciondelasetiquetasprimera}) over the set of $m-1$ elements $ \{ d_{1}, \sdots , d_{m} \} \setminus \{ d_{k'} \} $. This shows that we can omit repeated consecutive elements from the list $ d_{1}, \sdots , d_{m} $ and still get the same result. Hence, without loss of generality, we can assume that $d_{1}, \sdots , d_{m}$ are all different from their neighbors modulo $m$.

We proceed by induction over $m$. If $m=2$ then $$ \frac{(1-d_{1}d_{2})^{2}}{(1-d_{1}^{2})(1-d_{2}^{2})} = \frac{(1-d_{1}d_{2})^{2}}{(1-d_{1}d_{2})^{2}-(d_{1}-d_{2})^{2}} > 1 \ .$$\par
Now assume the statement holds for $ m-1$. First, we make an observation: There exists $k \in \{ 1, \sdots , m \}$ such that $d_{k}> d_{k-1} , d_{k+1}$ or $d_{k}< d_{k-1} , d_{k+1}$. To see this, suppose $d_{2} > d_{1}$. If $d_{3} < d_{2}$, then $k=2$ satisfies the property, and if not, then $d_{1}< d_{2}< d_{3}$. Similarly if $d_{4} < d_{3}$ we have finished, otherwise $d_{1}< d_{2}< d_{3} < d_{4}$. In the worst scenario $d_{1} < \scdots < d_{m}$, which implies that $d_{m} > d_{m-1},d_{1}$. The reasoning for when $d_{2} < d_{1}$ is analogous.\par
Let $k$ be such an index and define the quadratic polynomial 
$ p(x)=(x-d_{k-1})(x-d_{k+1})$. Notice that for $|x|\neq 1$ 
\begin{equation}\label{afirmaciondelasetiquetassegunda}
p(x) > 0 \ \Leftrightarrow \ \frac{(1-d_{k-1}x)(1-d_{k+1}x)}{1-x^{2}} > 1-d_{k-1}d_{k+1} \ . \nonumber
\end{equation}

The roots of $p$ are $d_{k-1}$ and $d_{k+1}$, and so $p(x)> 0$ for all $x > d_{k-1}, d_{k+1}$ or $x < d_{k-1} , d_{k+1}$. In particular $p(d_{k}) > 0$, therefore
$$ \frac{(1-d_{k-1}d_{k})(1-d_{k+1}d_{k})}{1-d_{k}^{2}} > 1-d_{k-1}d_{k+1} \ , $$
and $$ \frac{1-d_{k-1}d_{k}}{1-d_{k-1}^{2}} \cdot \frac{1-d_{k}d_{k+1}}{1-d_{k}^{2}} > \frac{1-d_{k-1}d_{k+1}}{1-d_{k-1}^{2}} \ . $$
Then
\begin{eqnarray}
\prod_{j=1}^{m} \frac{1-d_{j}d_{j+1}}{1-d_{j}^{2}} & > & \frac{1-d_{k-1}d_{k+1}}{1-d_{k-1}^{2}} \prod_{\substack{j=1 \\ j \neq k-1,k}}^{m} \frac{1-d_{j}d_{j+1}}{1-d_{j}^{2}}  \ , \nonumber
\end{eqnarray}
which by the induction hypothesis is greater than $1$ because, as before, the rightmost term is a product of the form of (\ref{afirmaciondelasetiquetasprimera}) over a set of $m-1$ elements.
\end{proof}
\end{description}
\begin{proof}[Proof of Lemma \ref{lemmadeladesigualdaddelafase}] 
For each $j = 1, \sdots , m$ we have
\begin{equation}\label{proposiciondespuesdelasetiquetascontrianguloycauchyschwarz}
| 1-\langle \xi_{j} , \xi_{j+1} \rangle | \geq 1-|\langle \xi_{j} , \xi_{j+1} \rangle| \geq 1-|\xi_{j}||\xi_{j+1}| \ ,
\end{equation}
which we use for
\begin{eqnarray}
\mathrm{Im}\big( \Phi(\xi_{1}, \sdots , \xi_{m}) \big) & = & \mathrm{Im}\left( i \Log \prod_{j=1}^{m} \frac{1-\langle \xi_{j} , \xi_{j+1} \rangle}{1-|\xi_{j}|^{2}} \right) \nonumber \\
& \ & \ \nonumber \\
& = & \ln \left\vert \prod_{j=1}^{m} \frac{1-\langle \xi_{j} , \xi_{j+1} \rangle}{1-|\xi_{j}|^{2}} \right\vert \nonumber \\
& \ & \ \nonumber \\
& \ge & \ln \prod_{j=1}^{m} \frac{1-|\xi_{j}|| \xi_{j+1}|}{1-|\xi_{j}|^{2}} \nonumber \ .
\end{eqnarray}
The last product is nonnegative by (\ref{afirmaciondelasetiquetasprimera}), thus proving the inequality part. The if part of the equality characterization is straightforward. On the other hand, if the equality holds, then by Claim \ref{proposiciondelasetiquetas} $ |\xi_{1}|=\scdots = |\xi_{m}| $. Inequality (\ref{proposiciondespuesdelasetiquetascontrianguloycauchyschwarz}) involves the Cauchy-Schwarz inequality, so the equality being attained implies that $ \xi_{j} $ and $\xi_{j+1}$ are linearly dependent, and thus $\xi_{j} = \pm \xi_{j+1}$. Finally, if 
$\xi_{j} =- \xi_{j+1}$, then
$$ | 1-\langle \xi_{j} , \xi_{j+1} \rangle | = 1+| \xi_{j} |^{2} > 1-|\xi_{j}||\xi_{j+1}| \ ,$$ which contradicts the assumption of equality, hence $\xi_{j}=\xi_{j+1}$ and $ \xi_{1}=\scdots = \xi_{m} $.
\end{proof}

\subsection{Asymptotic of the trace of a multiple composition}\label{laseccionconelcomportamientoasintotico}

This section is devoted to finding the asymptotic behavior of the integral expression in (\ref{latrazadelacomposicion}) when $\alpha$ tends to infinity. Throughout the section, we will keep $ m \geq 2 $ and $ a_{1}, \sdots , a_{m} \in C^{\infty}_{0}(\Gamma) $ fixed and denote $S:=T_{a_{1}\mathrm{d}\sigma} \circ \scdots \circ T_{a_{m}\mathrm{d}\sigma} $.\par

Consider $\gamma : B(0,1) \subset \mathbb{R}^{d} \rightarrow \Gamma$, \ $ t=(t_{1}, \sdots , t_{d})  \mapsto \gamma(t) $ a parametrization of an open set in $ \Gamma$. Let $t,s_{1}, \sdots , s_{m-1}$, $s_{q}=(s^{(q)}_{1}, \sdots , s^{(q)}_{d})$ denote $m$ variables in $\mathbb{R}^{d}$ and $ \pmb{s}=(s_{1}, \sdots , s_{m-1}) \in \mathbb{R}^{d(m-1)} $. We define
\begin{eqnarray}
\Xi & : & B(0,\varepsilon)^{m} \subset \mathbb{R}^{dm} \rightarrow \Gamma^{m} \nonumber \\
\Xi(t,\pmb{s}) & = & \big( \gamma(t), \gamma(t+s_{1}), \sdots , \gamma(t+s_{m-1}) \big) \ . \nonumber
\end{eqnarray}
Choosing $ \varepsilon $ small enough $\Xi$ is a parametrization for an open set $\mathcal{V} \subseteq \Gamma^{m}$, which intersects the diagonal of $ \Gamma^{m} $ whenever $ \pmb{s}=\pmb{0} $, i.e., $ s_{1} = \scdots = s_{m-1} = 0 $. Let $ \chi \in C_{0}^{\infty}(\mathcal{V}) $. We wish to compute

\begin{multline}
I_{\chi}(\alpha) = c_{\alpha}^{m} \int_{\mathcal{V}} e^{i\alpha \Phi(\xi_{1}, \dots , \xi_{m})} \left( \prod_{j=1}^{m} \frac{a_{j}(\xi_{j})}{(1-\langle \xi_{j} , \xi_{j+1} \rangle)^{n+1}} \right) \chi(\xi_{1}, \sdots , \xi_{m}) \\
\cdot \mathrm{d}\sigma(\xi_{1}) \cdots \mathrm{d}\sigma(\xi_{m}) \ ,
\end{multline}
asymptotically as $\alpha \rightarrow \infty$ using stationary phase method through the parametrization $ \Xi (t,\pmb{s})$. The stationary phase will be made over $ \pmb{s} $ with $t$ fixed and then we will integrate over $t$. Regarding $\mathcal{V}$ as a subset of the product manifold $\Gamma^{m}$, its Riemannian metric $ G_{\Xi} $ satisfies $ \det G_{\Xi}(t,\pmb{s}) = {\det G(t)} \cdot {\det G(t+s_{1})} \cdots {\det G(t+s_{m-1})} $. Make $\varphi (t,\pmb{s}) = \Phi \circ \Xi = \Phi \big( \gamma(t), \gamma(t+s_{1}), \sdots , \gamma(t+s_{m-1}) \big) $ so that we have

\begin{align}\label{unprimereqnarrayparausarfaseestacionaria}
    \frac{I_{\chi}(\alpha)}{c_{\alpha}^{m}} & = \int\limits_{B(0,\varepsilon)^{m}} e^{i\alpha \Phi \big( \Xi(t, \pmb{s})\big)} \left( \prod_{j=1}^{m} \frac{a_{j}(\gamma(t+s_{j-1}))}{\big( 1-\left\langle \gamma(t+s_{j-1}) , \gamma(t+s_{j}) \right\rangle \big)^{n+1}} \right) \chi \big( \Xi(t,\pmb{s}) \big) \nonumber \\
    & \quad \qquad \cdot \det G_{\Xi}(t,\pmb{s})^{\frac{1}{2}} \ \mathrm{d}(t,\pmb{s}) \nonumber \\
    & = \int\limits_{B(0,\varepsilon)} \left[ \int\limits_{B(0,\varepsilon)^{m-1}} e^{i\alpha \varphi(t, \pmb{s})} \left( \prod_{j=1}^{m} \frac{a_{j}(\gamma(t+s_{j-1}))\det G(t+s_{j})^{\frac{1}{2}}}{\big(1-\left\langle \gamma(t+s_{j-1}) , \gamma(t+s_{j}) \right\rangle \big)^{n+1}} \right) \right. \nonumber \\
    & \quad \qquad \cdot \chi \big( \Xi(t,\pmb{s})\big) \ \mathrm{d}\pmb{s} \left] \mathrm{d}t \phantom{ \int\limits_{B(0,\varepsilon)^{m-1}} } \right.
\end{align}

(For simplicity of the expression, we have adopted the convention that $ s_{0} := 0 =:s_{m} $). We claim that $ \varphi(t, \cdot) $ has a critic point at $\pmb{s}=\pmb{0}$. Let $ 1 \leq q \leq m-1 $ and $ 1 \leq j \leq d $. We differentiate $ \varphi(t, \cdot) $ with respect to $s^{(q)}_{j}$ using (\ref{laformadelaphifaseenterminosdellogdelak}) and that $\Log K(z,w)$ is holomorphic in $z$ and antiholomorphic in $w$. We obtain

{\allowdisplaybreaks
\begin{align}\label{laderivadadelaphiNUEVA}
    -i\frac{\partial}{\partial s^{(q)}_{j}} \varphi (t, \pmb{s}) & = \frac{1}{n+1} \sum_{\ell = 1}^{n} \Bigg[ \ \frac{\partial}{\partial z_{\ell}} \Log K\big(\gamma(t+s_{q}),\gamma(t+s_{q})\big) \cdot \partial_{j}\gamma_{\ell}(t+s_{q}) \nonumber \\
    & \quad \qquad + \ \frac{\partial}{\partial \overline{w}_{\ell}} \Log K\big(\gamma(t+s_{q}),\gamma(t+s_{q})\big) \cdot \partial_{j}\overline{\gamma_{\ell}}(t+s_{q}) \nonumber \\
    & \quad \qquad - \ \frac{\partial}{\partial z_{\ell}} \Log K\big(\gamma(t+s_{q}),\gamma(t+s_{q+1})\big) \cdot \partial_{j}\gamma_{\ell}(t+s_{q}) \nonumber \\
    & \quad \qquad - \ \frac{\partial}{\partial \overline{w}_{\ell}} \Log K\big(\gamma(t+s_{q-1}),\gamma(t+s_{q})\big) \cdot \partial_{j}\overline{\gamma_{\ell}}(t+s_{q}) \ \Bigg]
\end{align}
}

from which it follows that $ \nabla_{\pmb{s}}\varphi (t,\mathbf{0}) = 0 $ . \\
We now use (\ref{laderivadadelaphiNUEVA}) to compute the Hessian of $\varphi$ at $\pmb{s}=\pmb{0}$. Let $ 1 \leq q' \leq m-1 $ and $ 1 \leq k \leq d $. Notice that if $ q' \notin \{ q-1,q,q+1 \} $, then (\ref{laderivadadelaphiNUEVA}) does not depend on $ s^{(q')}_{k} $, and therefore $ \frac{\partial^{2}}{\partial s^{(q')}_{k} \partial s^{(q)}_{j}} \varphi (t, \pmb{s}) = 0 $. For the remaining three no null cases we have

\begin{equation}
    \frac{\partial^{2}}{\partial s^{(q-1)}_{k} \partial s^{(q)}_{j}} \varphi (t,\mathbf{0}) = \frac{-i}{n+1} \sum_{\ell , r = 1}^{n} \Bigg[ \ \frac{\partial^{2}}{\partial z_{r} \partial \overline{w}_{\ell}} \Log K\big(\gamma(t),\gamma(t)\big) 
     \partial_{j}\overline{\gamma_{\ell}}(t) \partial_{k}\gamma_{r}(t) \Bigg] \ , \nonumber
\end{equation}
\begin{equation}
    \frac{\partial^{2}}{\partial s^{(q+1)}_{k} \partial s^{(q)}_{j}} \varphi (t,\mathbf{0}) = \frac{-i}{n+1} \sum_{\ell , r = 1}^{n} \Bigg[ \ \frac{\partial^{2}}{\partial z_{\ell} \partial \overline{w}_{r}} \Log K\big(\gamma(t),\gamma(t)\big) 
    \partial_{j}\gamma_{\ell}(t) \partial_{k}\overline{\gamma_{r}}(t) \Bigg] \ , \nonumber
\end{equation}
\begin{equation}
\begin{split}
    i \frac{\partial^{2}}{\partial s^{(q)}_{k} \partial s^{(q)}_{j}} \varphi (t,\mathbf{0}) & = \frac{1}{n+1} \sum_{\ell , r = 1}^{n} \Bigg[ \ \frac{\partial^{2}}{\partial z_{\ell} \partial \overline{w}_{r}} \Log K\big(\gamma(t),\gamma(t)\big) \partial_{j}\gamma_{\ell}(t)\partial_{k}\overline{\gamma_{r}} \\
    & \quad \qquad + \ \frac{\partial^{2}}{\partial z_{r} \partial \overline{w}_{\ell}} \Log K\big(\gamma(t),\gamma(t)\big) \partial_{j}\overline{\gamma_{\ell}}(t)\partial_{k}\gamma_{r} \ \Bigg] \ .
\end{split} \nonumber
\end{equation}

Recall that $$ \frac{1}{n+1} \frac{\partial^{2}}{\partial z_{\ell} \partial \overline{z}_{r}} \log K(z,z)\Big\rvert_{z=\gamma(t)} $$
are the components $b_{\ell r}$ of the Bergman Hermitian metric of $\mathbb{B}_{n}$ at $\gamma(t)$. Comparing with (\ref{relaciondelasmetricasengamma}) and (\ref{relaciondelasmetricasengammaconjugada}) of Lemma \ref{lemaconlasrealcionesdelasformasbilinealesengamma} we get
{\allowdisplaybreaks
\begin{eqnarray}
-i \frac{\partial^{2}}{\partial s^{(q-1)}_{k} \partial s^{(q)}_{j}} \varphi (t,\mathbf{0}) & = & -g_{jk}(p) + i h_{jk}(p) \ , \nonumber \\
& \ & \ \nonumber \\
-i \frac{\partial^{2}}{\partial s^{(q+1)}_{k} \partial s^{(q)}_{j}} \varphi (t,\mathbf{0}) & = & -g_{jk}(p) - i h_{jk}(p) \ , \nonumber \\
& \ & \ \nonumber \\
-i \frac{\partial^{2}}{\partial s^{(q)}_{k} \partial s^{(q)}_{j}} \varphi (t,\mathbf{0}) & = & 2g_{jk}(p) \ , \nonumber
\end{eqnarray}
}
where in accordance with Section \ref{lasecciondelametricariemannianadeGamma}, $g_{jk}(p)$ and $h_{jk}(p)$ are the components at $p=\gamma(t)$ of the Rieamannian metric $g$ and the skew-symmetric form $h$ in $\Gamma$. If we define the matrices $ G=( g_{jk}(p)) $ and $ H=(h_{jk}(p)) $, then the Hessian matrix of $\varphi$ at $\pmb{s}=\pmb{0}$ can be expressed as a $(m-1)\times(m-1)$ block matrix with blocks of size $d\times d$ in the form
\small
\begin{equation}\label{formamatricialenbolquesdelhessianodelaphi}
-i \mathrm{Hess}\varphi(t,\mathbf{0}) = \begin{pmatrix}
2G & -G-i H & 0 & \scdots & 0 \\
-G+i H & 2G & -G-i H & \scdots & 0 \\
0 & -G+i H & 2G & \scdots & 0 \\
\vdots & \vdots & \vdots & \vdots & \vdots \\
0 & \scdots & -G+i H & 2G & -G-i H \\
0 & \scdots & 0 & -G+i H & 2G
\end{pmatrix}  . \nonumber
\end{equation}
\normalsize

To express the determinant of the Hessian, make $W:=G^{-1}H$ to rewrite

\small
\begin{equation}\label{formamatricialdelquessianodelaphifactorizada}
{ -i \mathrm{Hess}\varphi(t,\mathbf{0}) {\scriptstyle =} \begin{pmatrix}
G & 0 & {\scriptstyle\cdots} & 0 \\
0 & G & {\scriptstyle\cdots} & 0 \\
{\scriptstyle\vdots} & {\scriptstyle\vdots} & {\scriptstyle\vdots} & {\scriptstyle\vdots} \\
0 & {\scriptstyle\cdots} & G & 0 \\
0 & {\scriptstyle\cdots} & 0 & G 
\end{pmatrix}
\begin{pmatrix}
2I & -I-i W & {\scriptstyle\cdots} & 0 \\
-I+i W & 2I  & {\scriptstyle\cdots} & 0 \\
{\scriptstyle\vdots} & {\scriptstyle\vdots} & {\scriptstyle\vdots} & {\scriptstyle\vdots} \\
0 & {\scriptstyle\cdots} & 2I & -I-i W \\
0 & {\scriptstyle\cdots} & -I+i W & 2I 
\end{pmatrix} , }
\end{equation}
\normalsize

and let $\mathcal{O}_{m-1}$ be the rightmost matrix of the factorization. We then have
\begin{equation}
\det \left( -i\mathrm{Hess}\varphi(t,\mathbf{0}) \right) =( \det G )^{m-1} \det \mathcal{O}_{m-1} \ . \nonumber
\end{equation}

In Appendix \ref{laseccionconloscalculosparaeldeterminante} it is shown that $\det \mathcal{O}_{m-1}$ is positive, while so is $\det G$ because $G$ is positive-definite. This implies that $\det \left( -i\mathrm{Hess}\varphi(t,\mathbf{0}) \right)$ is also positive and so the application of the stationary phase method is justified. Apply Theorem 7.7.5 \cite{hormander} to estimate the integral inside the brackets in (\ref{unprimereqnarrayparausarfaseestacionaria}), to obtain

\begin{equation}
\begin{split}
    & \int\limits_{B(0,\varepsilon)^{m-1}} e^{i\alpha \varphi(t, \pmb{s})} \left( \prod_{j=1}^{m} \frac{a_{j}(\gamma(t+s_{j-1}))\det G(t+s_{j})^{\frac{1}{2}}}{\big(1-\left\langle \gamma(t+s_{j-1}) , \gamma(t+s_{j}) \right\rangle \big)^{n+1}} \right) \chi \big( \Xi(t,\pmb{s})\big) \ \mathrm{d}\pmb{s} = \\
    & \quad \qquad \mathrm{det}\left( \frac{\alpha}{2\pi i} \mathrm{Hess}\varphi(t, \pmb{0}) \right)^{-\frac{1}{2}} \left[ \ \prod_{j=1}^{m} \frac{a_{j}(\gamma(t))\det G(t)^{\frac{1}{2}}}{\big(1-| \gamma(t) |^{2} \big)^{n+1}} \ \right] \ \chi\big( \gamma(t), \sdots , \gamma(t)\big) \\
    & \quad \qquad \qquad \cdot \left( 1 + O\left( \alpha^{-1} \right) \right) \ .
\end{split} \nonumber
\end{equation}

Notice that up until this point we had not made use of the smoothness of the symbols $a_{1}, \sdots , a_{m}$. In order to apply the cited theorem to obtain the previous estimate, we need each $ a_{j} \in C^{2k}_{0}(\Gamma) $ for some $k \geq \frac{d(m-1)}{2} +1 $. As we wish for the hypotheses on the symbols to be independent of their quantity $m$, we have rather requested $a_{j} \in C^{\infty}_{0}(\Gamma) $ for each $j$. Also, as we will see later, this same $m$ will represent the exponent of powers of $\topnorm$. This is the reason why we requested $a \in C^{\infty}_{0}(\Gamma) $ in Definition \ref{definiciontoeplitz}.

Substituting on (\ref{unprimereqnarrayparausarfaseestacionaria}) the resulting estimate we get

\begin{align}\label{uneqnarraylargoenelcalculoasintotico}
    \frac{I_{\chi}(\alpha)}{c_{\alpha}^{m}} & = \int\limits_{B(0,\varepsilon)} \Bigg[ \det\left( \frac{\alpha}{2\pi i} \mathrm{Hess}\varphi(t,0) \right)^{-\frac{1}{2}} \left( \prod_{j=1}^{m} \frac{a_{j}(\gamma(t))\det G(t)^{\frac{1}{2}}}{(1-\left| \gamma(t) \right|^{2} )^{n+1}} \right) \nonumber \\
    & \quad \qquad \cdot \chi( \gamma(t), \sdots , \gamma(t)) \left( 1 + O\left(\alpha^{-1}\right) \right) \ \Bigg] \mathrm{d}t \nonumber \\
    & \ \nonumber \\
    & = \left(\frac{2\pi}{\alpha}\right)^{\frac{d(m-1)}{2}} \int\limits_{B(0,\varepsilon )} \left( \prod_{j=1}^{m} a_{j}(\gamma(t)) \right) \cdot \frac{\chi ( \gamma(t), \sdots , \gamma(t))}{\left( 1-|\gamma(t)|^{2} \right)^{m(n+1)} } \nonumber \\
    & \quad \qquad \cdot \frac{\det G(t)^{\frac{1}{2}} }{(\det \mathcal{O}_{m-1})^{\frac{1}{2}}} \mathrm{d}t \left( 1 + O\left( \alpha^{-1} \right) \right) \nonumber \\
    & \ \nonumber \\
    & = \left(\frac{2\pi}{\alpha}\right)^{\frac{d(m-1)}{2}} \int\limits_{\gamma (B(0,\varepsilon))} \left( \prod_{j=1}^{m} a_{j}(\xi) \right) \frac{\chi (\xi, \sdots , \xi)}{\left( 1-|\xi|^{2} \right)^{m(n+1)} } \frac{\mathrm{d}\sigma(\xi)}{(\det \mathcal{O}_{m-1})^{\frac{1}{2}}} \nonumber \\
    & \quad \qquad \cdot \left(1 + O\left( \alpha^{-1} \right) \right) \ .
\end{align}

In this last expression, the integration can be changed to be done on all $\Gamma$ instead of just $\gamma (B(0,\varepsilon))$ without affecting the outcome. To see this, note that if $\xi \notin \gamma(B(0, \varepsilon))$ then $ (\xi, \sdots , \xi) \notin \Xi(B(0, \varepsilon)^{m}) = \mathcal{V} \supset \mathrm{supp} \chi $.

Let $ \gamma_{1}, \sdots , \gamma_{L}: B(0,1)\subset \mathbb{R}^{d} \rightarrow \Gamma $ be parametrizations of $\Gamma$ such that $\cup_{\ell=1}^{L}\gamma_{\ell}( B(0,1) ) = \Gamma$. 
Without loss of generality, we have taken this as a finite union since, even if $\Gamma$ is not compact, the integral in which we are interested is done on the supports of the functions $a_{j}$. Take $\varepsilon$ small enough and define
$$ \Xi_{\ell} : B(0,\varepsilon)^{m} \rightarrow \Gamma^{m} $$ $$\Xi_{\ell}(t,\pmb{s}) = \big( \gamma_{\ell}(t), \gamma_{\ell}(t+s_{1}), \sdots , \gamma_{\ell}(t+s_{m-1}) \big) \ , $$
so that each $ \Xi_{\ell}( B(0,\varepsilon)^{m} )$ is, as described in the previous construction, an open set in $\Gamma^{m}$ that intersects the diagonal when $\pmb{s}=\mathbf{0}$. Then 
$$ \mathcal{U}:= \bigcup_{\ell=1}^{L} \Xi_{\ell}( B(0,\varepsilon)^{m} ) $$ is an open neighborhood of the diagonal in $\Gamma^{m}$ and $ \left\{ \Xi_{\ell} , B(0,\varepsilon )^{m} \right\}_{\ell=1, \dots , L} $ is an atlas for $\mathcal{U}$. It follows from (\ref{formulatraza}) and Lemma \ref{lemmadeladesigualdaddelafase} that $\mathrm{Tr}(S)$ decays exponentially fast outside $\mathcal{U}$ when $\alpha $ grows. Let $\{ \chi_{\ell} \}_{\ell=1, \dots , L}$ be a partition of unity subordinated to the aforementioned atlas. Applying what we obtained in (\ref{uneqnarraylargoenelcalculoasintotico}) to each $\chi_{\ell}$ and substituting it into (\ref{formulatraza}) we get

{\allowdisplaybreaks
\begin{align}\label{asintotadelatraza}
    \mathrm{Tr}(S) & = c_{\alpha}^{m} \int\limits_{ \mathcal{U}} e^{i\alpha \Phi(\xi_{1}, \dots , \xi_{m})} { \left( \prod_{j=1}^{m} \frac{a_{j}(\xi_{j})}{(1-\langle \xi_{j} , \xi_{j+1} \rangle)^{n+1}} \right)} \ \mathrm{d}\sigma(\xi_{1}) \cdots \mathrm{d}\sigma(\xi_{m}) \nonumber \\
    & \quad \qquad + O\left( \alpha^{-\infty} \right) \nonumber \\
    & = c_{\alpha}^{m} \sum_{\ell=1}^{L} \int\limits_{ \Xi_{\ell}( B(0,\varepsilon)^{ m} )} e^{i\alpha \Phi(\xi_{1}, \dots , \xi_{m})} { \left( \prod_{j=1}^{m} \frac{a_{j}(\xi_{j})}{(1-\langle \xi_{j} , \xi_{j+1} \rangle)^{n+1}} \right)} \chi_{\ell}({ \xi_{1}, {\scriptstyle \dots } , \xi_{m}}) \nonumber \\
    & \quad \qquad \cdot \mathrm{d}\sigma {\scriptstyle(\xi_{1})} \cdots \mathrm{d}\sigma{\scriptstyle (\xi_{m})} + O\left( \alpha^{-\infty} \right) \nonumber \\
    & = \sum_{\ell=1}^{L} \ I_{\chi_{\ell}}(\alpha) + O\left( \alpha^{-\infty} \right) \nonumber \\
    & = c_{\alpha}^{m} \left(\frac{2\pi}{\alpha}\right)^{\frac{d(m-1)}{2}} \sum_{\ell=1}^{L} \int\limits_{\Gamma} \left( \prod_{j=1}^{m} a_{j}(\xi) \right) \frac{\chi_{\ell}(\xi, \sdots , \xi)}{\left( 1-|\xi|^{2} \right)^{m(n+1)} } \frac{\mathrm{d}\sigma(\xi)}{(\det \mathcal{O}_{m-1})^{\frac{1}{2}}} \nonumber \\
    & \quad \qquad \cdot \left( 1 + O\left(\alpha^{-1}\right) \right) \nonumber \\
    & = c_{\alpha}^{m} \left(\frac{2\pi}{\alpha}\right)^{\frac{d(m-1)}{2}} \int\limits_{\Gamma} \left( \prod_{j=1}^{m} a_{j}(\xi) \right) \frac{\sum_{\ell=1}^{L} \chi_{\ell}(\xi, \sdots , \xi)}{\left( 1-|\xi|^{2} \right)^{m(n+1)} } \frac{\mathrm{d}\sigma(\xi)}{(\det \mathcal{O}_{m-1})^{\frac{1}{2}}} \nonumber \\
    & \quad \qquad \cdot \left( 1 + O\left(\alpha^{-1}\right) \right) \nonumber \\
    & = c_{\alpha}^{m} \left(\frac{2\pi}{\alpha}\right)^{\frac{d(m-1)}{2}} \int\limits_{\Gamma} \left( \prod_{j=1}^{m} \frac{a_{j}(\xi)}{\left( 1-|\xi|^{2} \right)^{n+1}} \right) \frac{\mathrm{d}\sigma(\xi)}{(\det \mathcal{O}_{m-1})^{\frac{1}{2}}} \left( 1 + O\left(\alpha^{-1}\right) \right) ,
\end{align}
}
since $ \sum_{\ell} \chi_{\ell} = 1 $  in $\Gamma$.

In order to obtain an explicit formula from (\ref{asintotadelatraza}), we would like $\det \mathcal{O}_{m-1}$ to be computable, while it would be ideal for it to be constant for all $\xi = \gamma(t) \in \Gamma$. This does not need to be true in general, but in Appendix \ref{laseccionconloscalculosparaeldeterminante} is shown this is the case when $\Gamma$ is isotropic or co-isotropic. Namely, Proposition \ref{elcorolariodondeconlcuimoseldeterminantedelaoexplicito} states
\begin{equation}
\sqrt{ \det \mathcal{O}_{m-1} } \ = \ \left\{ \begin{array}{cl}
m^{\frac{d}{2}} & \Gamma \mbox{ isotropic} \\
2^{(d-n)(m-1)}m^{n-\frac{d}{2}} & \Gamma \mbox{ co-isotropic}
\end{array} \right. \ ,
\end{equation}
which, rewriting it in terms of $d'$, leads to
\begin{equation}\label{launificaciondelosdeterminantesconladapostrofe}
2^{\frac{d(m-1)}{2}}( \mathrm{det} \mathcal{O}_{m-1} )^{-\frac{1}{2}} =  2^{\frac{d'(m-1)}{2}} m^{-\frac{d'}{2}}
\end{equation}
for both $\Gamma$ isotropic or co-isotropic. Also, recall that $ \topnorm = \frac{n!}{2^{\frac{d'}{2}}\pi^{\frac{d}{2}}}\alpha^{-n+\frac{d}{2}} T_{a\mathrm{d}\sigma} $. We conclude this section by presenting a statement that describes the resulting estimation when adding this hypothesis.

\begin{lemma}{\bf [Asymptotics for the trace of a composition]}\label{lapropodelaasintoticadelatrazadelacomposicion}
Let $\Gamma$ be isotropic or co-isotropic and $a_{1}, \sdots , a_{m} \in C^{\infty}_{0}(\Gamma)$. Then
\begin{equation}
     \mathrm{Tr}\left( \widehat{T}_{a_{1}\mathrm{d}\sigma} \circ \scdots \circ \widehat{T}_{a_{m}\mathrm{d}\sigma} \right) = \frac{\left( \frac{\alpha}{\pi} \right)^{\frac{d}{2}}}{(2m)^{\frac{d'}{2}}} \int_{\Gamma} \prod_{j=1}^{m} \frac{a_{j}(\xi)}{(1-|\xi|^{2})^{n+1}} \ \mathrm{d}\sigma(\xi) 
     \left( 1+O\left( \frac{1}{\alpha} \right) \right) . \nonumber
\end{equation}
\end{lemma}
\begin{proof}
It can be verified using Stirling's formula that $$ c_{\alpha} = \frac{\Gamma(\alpha+1+n)}{n!\Gamma(\alpha+1)} = \frac{\alpha^{n}}{n!}\left( 1+O\left( \frac{1}{\alpha} \right) \right) . $$
By the linearity of the trace, the result will follow from (\ref{asintotadelatraza}) by substituting the last identity for $c_{\alpha}$ and the expression for $ 2^{\frac{d(m-1)}{2}}( \mathrm{det} \mathcal{O}_{m-1} )^{-\frac{1}{2}} $ given in (\ref{launificaciondelosdeterminantesconladapostrofe}).
\end{proof}

\section{Proof of the main result}
\label{laseccionconelteorema}

In this section, we demonstrate Theorem \ref{teoremaszego}. First, we state the result for monomial functions and then extend it to polynomials by linearity. Next, we will use a Weierstrass argument to derive the theorem for continuous functions.\par

\subsection{The polynomial case}

As a consequence of Lemma \ref{lapropodelaasintoticadelatrazadelacomposicion} we obtain the following proposition, which can be seen as a first form of our Szeg\"{o} limit theorem applying to monomial functions $\varphi(t)=t^{m}$.

\begin{proposition}\label{corolariomonomio}
Let $\Gamma$ be isotropic or co-isotropic and $a \in C^{\infty}_{0}(\Gamma)$. Then
\begin{equation}\label{ecuacionencorolariomonomio}
\mathrm{Tr}\left( \left(\widehat{T}_{a\mathrm{d}\sigma}\right)^{m} \right) = \frac{\left( \frac{\alpha}{\pi} \right)^{\frac{d}{2}}}{(2m)^{\frac{d'}{2}}} \int_{\Gamma} \left( \frac{a(\xi)}{(1-|\xi|^{2})^{n+1}} \right)^{m} \ \mathrm{d}\sigma(\xi) \ \left( 1+O\left( \frac{1}{\alpha} \right) \right) \ .
\end{equation}
\end{proposition}
\begin{proof}
Take $a_{1}= \sdots = a_{m} = a$ in Lemma \ref{lapropodelaasintoticadelatrazadelacomposicion}.
\end{proof}

\subsubsection{An estimate for the operator norm}

As a first application, we obtain an estimation of the growth rate of the operator norm of the Toeplitz operators. This result shows that the operators defined in (\ref{lanormalizacionprimeravez}) are indeed normalizations of the Toeplitz operators defined in Definition \ref{definiciontoeplitz}.

\begin{corollary}\label{corolariodelanormalizacion}
Let $\Gamma$ be isotropic or co-isotropic and $a \in C^{\infty}_{0}(\Gamma)$. There exists a constant $C>0$ such that
\begin{equation}
\big|\big| \widehat{T}_{a\mathrm{d}\sigma} \big|\big|_{\scriptscriptstyle \mathcal{L}( A^{2}_{\alpha}(\mathbb{B}_{n}) )} < C \nonumber
\end{equation}
for all $\alpha > 0$.
\end{corollary}
\begin{proof}
Let $a=a_{1}+i a_{2}$ with $a_{1}$ and $a_{2}$ real valued. We can write
$$ \widehat{T}_{a \mathrm{d}\sigma}  \ = \ \widehat{T}_{(a_{1}+||a||_{\infty}) \mathrm{d}\sigma} \ + \ i \widehat{T}_{(a_{2}+||a||_{\infty}) \mathrm{d}\sigma} \ - \ (1+i) \widehat{T}_{||a||_{\infty} \mathrm{d}\sigma} \ , $$
and so
\begin{equation}
\begin{split}
   \big| \big| \widehat{T}_{a \mathrm{d}\sigma} \big| \big|_{\scriptscriptstyle \mathcal{L}( A^{2}_{\alpha}(\mathbb{B}_{n}))} & \leq \big| \big| \widehat{T}_{ (a_{1}+||a||_{\infty}) \mathrm{d}\sigma} \big| \big|_{\scriptscriptstyle \mathcal{L}( A^{2}_{\alpha}(\mathbb{B}_{n}) )} + \big| \big| \widehat{T}_{ (a_{2}+||a||_{\infty}) \mathrm{d}\sigma} \big| \big|_{\scriptscriptstyle \mathcal{L}( A^{2}_{\alpha}(\mathbb{B}_{n}) )} \\
   & \quad \qquad + \sqrt{2} \big| \big| \widehat{T}_{ ||a||_{\infty} \mathrm{d}\sigma} \big| \big|_{\scriptscriptstyle \mathcal{L}( A^{2}_{\alpha}(\mathbb{B}_{n}) )} \ ,
\end{split} \nonumber
\end{equation}
where the three operators on the right-hand side have positive symbols. Therefore, it suffices to prove the proposition for $a\geq 0$, which in turn implies that $\widehat{T}_{a \mathrm{d}\sigma}$ is self-adjoint and positive. \par
Let $a \geq 0$ and $\ell$ be a positive integer.
Then
\begin{equation}\label{desigualdadnormasendemostraciondelanormalizacion}
\big| \big| \widehat{T}_{a \mathrm{d}\sigma} \big| \big|^{\ell}_{\scriptscriptstyle \mathcal{L}( A^{2}_{\alpha}(\mathbb{B}_{n}) )} = \big| \big| \big( \widehat{T}_{a \mathrm{d}\sigma} \big)^{\ell} \big| \big|_{\scriptscriptstyle \mathcal{L}( A^{2}_{\alpha}(\mathbb{B}_{n}) )} \leq \mathrm{Tr}\left( \big( \widehat{T}_{a \mathrm{d}\sigma} \big)^{\ell} \right) \ .
\end{equation}
On the other hand, using the result for the trace of the power of a Toeplitz operator in Proposition \ref{corolariomonomio},
\begin{eqnarray}
\mathrm{Tr}\left( \big( \widehat{T}_{a \mathrm{d}\sigma} \big)^{\ell} \right) & = & \frac{\left( \frac{\alpha}{\pi} \right)^{\frac{d}{2}}}{(2\ell)^{\frac{d'}{2}}} \int_{\Gamma} \left( \frac{a(\xi)}{(1-|\xi|^{2})^{n+1}} \right)^{\ell} \ \mathrm{d}\sigma(\xi) \ \left( 1+O\left( \frac{1}{\alpha} \right) \right) \nonumber \\
& \ & \ \nonumber \\
& \leq & C \left[ \sup\limits_{\xi \in \Gamma} \frac{a(\xi)}{(1-|\xi|^{2})^{n+1}} \right]^{\ell} \frac{\alpha^{\frac{d}{2}}}{\ell^{\frac{d'}{2}}} \ \nonumber ,
\end{eqnarray}
with the constant $C$ not depending on $\alpha$ nor $\ell$. Substitute this into (\ref{desigualdadnormasendemostraciondelanormalizacion}) and take $\ell$-th roots to get

\begin{equation}\label{ecuacioncondesigualdaddentrodelademodelanormalizacion}
\big| \big| \widehat{T}_{a \mathrm{d}\sigma} \big| \big|_{\scriptscriptstyle \mathcal{L}( A^{2}_{\alpha}(\mathbb{B}_{n}) )} \ \leq \ C^{\frac{1}{\ell}} \left[ \sup\limits_{\xi \in \Gamma} \frac{a(\xi)}{(1-|\xi|^{2})^{n+1}} \right] \frac{\alpha^{\frac{d}{2\ell}}}{\ell^{ \frac{d'}{2\ell}}} \ \leq \ C \left(\alpha^{\frac{1}{\ell}} \right)^{\frac{d}{2}} \ ,
\end{equation}
again with the last constant $C$ not depending on $\alpha$ nor $\ell$, and where we have used that $\left( \frac{1}{\ell}\right)^{\frac{1}{\ell}} \leq 1 $ for all $\ell$. This is valid for $\ell$ any positive integer, particularly for any $\ell > \alpha$. Then
$$ \big| \big| \widehat{T}_{a \mathrm{d}\sigma} \big| \big|_{\scriptscriptstyle \mathcal{L}( A^{2}_{\alpha}(\mathbb{B}_{n}) )} \ \leq \ C \inf\limits_{\ell \in \mathbb{N}} \left( \alpha^{\frac{1}{\ell}} \right)^{\frac{d}{2}} \ \leq \ C \left( \alpha^{\frac{1}{\alpha}} \right)^{\frac{d}{2}} \ \leq \ C \ . $$
Here we have used $ \lim\limits_{\alpha \rightarrow \infty} \alpha^{\frac{1}{\alpha}} =1 $. 
\end{proof}

\begin{remark}\label{remarkdespuesdelanormalizacion}
From (\ref{ecuacioncondesigualdaddentrodelademodelanormalizacion}) and the arguments given thereafter we find that if $a>0$, $\sup\limits_{\xi \in \Gamma} \frac{a(\xi)}{(1-|\xi|^{2})^{n+1}}$
is a bound for both the operator norm of $\widehat{T}_{a \mathrm{d}\sigma}$ and its largest eigenvalue.
\end{remark}

\subsubsection{The Szeg\"{o} limit theorem for polynomials}

We rely on the operator $\mathcal{Q}_{\varepsilon}$ to extend the result in Proposition \ref{corolariomonomio}, first to encompass polynomials. Recall that $\mathcal{Q}_{\varepsilon}$ is defined by
$$ \mathcal{Q}_{\varepsilon}(\varphi)(t) = \frac{1}{\Gamma(\varepsilon)} \int_{0}^{t} \varphi(s) \left( \ln \frac{t}{s} \right)^{\varepsilon-1} \ \frac{\mathrm{d}s}{s} $$
when $\varepsilon > 0$ and the identity when $\varepsilon = 0$.

\begin{myremark}
We summarize some important features of this operator.
\begin{itemize}
\item $\mathcal{Q}_{\varepsilon}$ is linear in the first entry.
\item For $p>0$, $ \mathcal{Q}_{\varepsilon}(s^{p})(t) = \frac{t^{p}}{p^{\varepsilon}} $. This follows from the change of variables ${x=p\ln\left( \frac{t}{s} \right)}$. In particular $ \mathcal{Q}_{\frac{d'}{2}}(s^{m})(t) = \frac{t^{m}}{m^{\frac{d'}{2}}} $. This remains valid even for $t\in\mathbb{C}$. In this case, the integration domain is taken to be the open line segment joining $0$ and $t$.
\item Equation (\ref{ecuacionencorolariomonomio}) in Proposition \ref{corolariomonomio} can be rewritten as 
\begin{multline}\label{ecuacionencorolariomonomioformaconq}
\mathrm{Tr}\left( \left(\widehat{T}_{a\mathrm{d}\sigma}\right)^{m} \right) = \frac{\left( \frac{\alpha}{\pi} \right)^{\frac{d}{2}}}{2^{\frac{d'}{2}}} \int_{\Gamma} \mathcal{Q}_{\frac{d'}{2}}\left(s^{m}\right) \left( a(\xi)(1-|\xi|^{2})^{-(n+1)} \right) \\ \mathrm{d}\sigma(\xi) \left( 1+O\left( \alpha^{-1} \right) \right) \ . 
\end{multline}
\end{itemize}
\end{myremark}

We now state the result for the polynomial case.

\begin{proposition}\label{corolariopolinomio}{\bf [Szeg\"{o} limit theorem for polynomials]}
Let $\Gamma$ be isotropic or co-isotropic, $a \in C^{\infty}_{0}(\Gamma)$, and $P$ be a polynomial without constant term. Then
\begin{equation}
\mathrm{Tr}\left( P\left(\widehat{T}_{a\mathrm{d}\sigma}\right) \right) = \frac{\left( \frac{\alpha}{\pi} \right)^{\frac{d}{2}}}{2^{\frac{d'}{2}}} \int_{\Gamma} \quessianobb{P} \ \mathrm{d}\sigma(\xi) \left( 1+O\left( \frac{1}{\alpha} \right) \right) \ . \nonumber
\end{equation}
\end{proposition}
\begin{proof}
The proof follows from (\ref{ecuacionencorolariomonomioformaconq}) and the linearity of both the trace and $\mathcal{Q}_{\varepsilon}$.
\end{proof}

\begin{remark}
Notice that for this case the hypothesis regarding $a>0$ in Theorem \ref{teoremaszego} is not needed and is therefore omitted.
\end{remark}

\subsection{The Szeg\"{o} limit theorem}\label{lasubsecciondelszegolimittheorem}

We now wish to extend our result to apply to more general continuous functions. From now on, we will consider $a>0$. Fix $R>0$ so that
$$ \sup\limits_{\xi \in \Gamma} \frac{a(\xi)}{(1-|\xi|^{2})^{n+1}} \leq R \ .$$
From Remark \ref{remarkdespuesdelanormalizacion} we infer that $$ \bigcup_{\alpha>0} \sigma\left( \topnorm \right) \subseteq [0,R] \ . $$
This is important as $\varphi(\topnorm)$ is only well defined if $\varphi$ is continuous on the spectrum of $\topnorm$. Because we seek to analyze the asymptotic behavior of $\varphi(\topnorm)$ when $\alpha$ grows, we require $\varphi$ to be continuous on the spectrum of $\topnorm$ for all $\alpha$.\par

We define the family of functions for which we will extend the past results.
\begin{definition}
$$ \mathcal{C} = \left\{ \varphi \in C[0,R] \ : \ \exists p>0 \ , \ \frac{\varphi(t)}{t^{p}} \in C[0,R] \right\} \ . $$
\end{definition}

\begin{myremark}
Without loss of generality, we can restrict the definition to consider only ${p \in (0,1)}$.
\end{myremark}

\begin{proposition}\label{propoconloscuatroincisosdelaclaseC}
The family $\mathcal{C}$ satisfies: 
\begin{itemize}
\item[a)] $\mathcal{C} \supset \{ P \ : \ \text{P is a polynomial without constant term} \}$ .
\item[b)] $ \varphi \left( \widehat{T}_{a\mathrm{d}\sigma} \right) $ is trace class for every $\varphi \in \mathcal{C}$.
\item[c)] $ \mathcal{Q}_{\varepsilon}(\varphi)(t) $ is well defined $\forall \ \varepsilon \geq 0$, $ t \in [0,R] $ and $ \varphi \in \mathcal{C}$. Moreover, if $p>0$ such that $ \frac{\varphi(s)}{s^{p}} \in C[0,R]$, then
\begin{equation}
|\mathcal{Q}_{\varepsilon}(\varphi)(t)| \leq \norma{\frac{\varphi}{s^{p}}}{L^{\infty}[0,R]} \frac{t^{p}}{p^{\varepsilon}} \ . \nonumber
\end{equation}
\item[d)] For every $\varphi \in \mathcal{C}$ $$ \int\limits_{\Gamma} \quessianobb{\varphi} \ \mathrm{d}\sigma(\xi) < \infty  \ .$$
\end{itemize}
\end{proposition}
\begin{proof}
{\it a)} is obvious. For {\it b)} let $ \lambda_{j} $ be the eigenvalues of $\widehat{T}_{a\mathrm{d}\sigma}$. Recall that as $a$ is taken to be positive, so is $\topnorm$. Then for $\varphi \in \mathcal{C}$ the norm of $ \varphi \left( \widehat{T}_{a\mathrm{d}\sigma} \right) $ in the trace class $S_{1}\left( A^{2}_{\alpha}(\mathbb{B}_{n}) \right)$ satisfies
{\allowdisplaybreaks
\begin{eqnarray}
\left\vert \left\vert \varphi \left( \widehat{T}_{a\mathrm{d}\sigma} \right) \right\vert \right\vert_{\scriptscriptstyle S_{1}( A^{2}_{\alpha}(\mathbb{B}_{n}) )} & = & \sum_{j=1}^{\infty} |\varphi(\lambda_{j})| \nonumber \\
& \ & \ \nonumber \\
& = & \sum_{j=1}^{\infty} \frac{|\varphi(\lambda_{j})|}{|\lambda_{j}|^{p}} |\lambda_{j}|^{p} \nonumber \\
& \ & \ \nonumber \\
& \leq & \sup\limits_{t \in [0,R]} \left\vert \frac{\varphi(t)}{t^{p}} \right\vert \sum_{j=1}^{\infty} |\lambda_{j}|^{p} \nonumber \\
& \ & \ \nonumber \\
& = & \norma{\frac{\varphi}{s^{p}}}{L^{\infty}[0,R]} \ \left\vert \left\vert \widehat{T}_{a\mathrm{d}\sigma} \right\vert \right\vert_{\scriptscriptstyle S_{p}( A^{2}_{\alpha}(\mathbb{B}_{n}) )}^{p} \ , \nonumber
\end{eqnarray}
}
which is finite because $\topnorm$ belongs to the Schatten class $S_{p}( A^{2}_{\alpha}(\mathbb{B}_{n}) )$.
For {\it c)} we compute
\begin{eqnarray}
\left\vert \mathcal{Q}_{\varepsilon}(\varphi)(t) \right\vert & = & \left\vert \frac{1}{\Gamma(\varepsilon)} \int_{0}^{t} \varphi(s) \left( \ln \frac{t}{s} \right)^{\varepsilon-1} \ \frac{\mathrm{d}s}{s} \right\vert \nonumber \\
& \ & \ \nonumber \\
& \leq & \frac{1}{\Gamma(\varepsilon)} \int_{0}^{t} \left\vert \frac{\varphi(s)}{s^{p}} \right\vert s^{p} \left( \ln \frac{t}{s} \right)^{\varepsilon-1} \ \frac{\mathrm{d}s}{s} \nonumber \\
& \ & \ \nonumber \\
& \leq & \frac{\norma{\frac{\varphi}{s^{p}}}{L^{\infty}[0,R]}}{\Gamma(\varepsilon)} \int_{0}^{t} s^{p} \left( \ln \frac{t}{s} \right)^{\varepsilon-1} \ \frac{\mathrm{d}s}{s} \ = \ \norma{\frac{\varphi}{s^{p}}}{L^{\infty}[0,R]} \frac{t^{p}}{p^{\varepsilon}} \nonumber \ .
\end{eqnarray}

Finally, for {\it d)} we use {\it c)} to get
\begin{align}\label{desigualdadparausarenlaqconcirculo}
    \left\vert \int\limits_{\Gamma} \quessianobb{\varphi} \mathrm{d}\sigma(\xi) \right\vert & \leq \int\limits_{\Gamma} \left\vert \quessianobb{\varphi} \right\vert \mathrm{d}\sigma(\xi) \nonumber \\
    \leq \norma{\frac{\varphi}{s^{p}}}{L^{\infty}[0,R]} & \frac{1}{p^{\frac{d'}{2}}} \int\limits_{\Gamma} \left( \frac{a(\xi)}{1-|\xi|^{n+1}} \right)^{p} \ \mathrm{d}\sigma(\xi) \nonumber \\
    \leq \norma{\frac{\varphi}{s^{p}}}{L^{\infty}[0,R]} &\frac{\sigma(\mathrm{supp} (a))}{p^{\frac{d'}{2}}} \left( \sup\limits_{\xi \in \Gamma} \frac{a(\xi)}{1-|\xi|^{n+1}}  \right)^{p} < \infty ,
\end{align}

where $ \sigma(\mathrm{supp} (a)) = \int_{\mathrm{supp}(a)}  \ \mathrm{d}\sigma(\xi) $.
\end{proof}

For the proof of the Szeg\"{o} limit theorem we will use the following two results. The proof of the first one relies strongly on the hyperbolic geometry of $\mathbb{B}_{n}$ and is left for Appendix \ref{seccioncondemodelaclaimdepchica}.
\begin{description}
\item \begin{lemma}\label{TRAZAPCHICA}
Let $\Gamma \subset \mathbb{B}_{n}$ be a submanifold and $a \in C^{\infty}_{0}(\Gamma)$ a positive function. For $p \in (0,1)$
\begin{equation}
\left(\frac{\pi}{\alpha}\right)^{\frac{d}{2}} \mathrm{Tr}\left( \left(\widehat{T}_{a\mathrm{d}\sigma}\right)^{p} \right) \nonumber
\end{equation}
is bounded as $\alpha \rightarrow \infty$. 
\end{lemma}
\item \begin{lemma}\label{claimproductoconpolinomio}
Let $\Gamma$ be isotropic or co-isotropic, $a \in C^{\infty}_{0}(\Gamma)$ a positive function, and $P$ a polynomial without constant term. For $p \in (0,1)$,
\begin{multline*}
    \lim\limits_{\alpha \rightarrow \infty} \left(\frac{\pi}{\alpha}\right)^{\frac{d}{2}} \mathrm{Tr}\left[ \left(\widehat{T}_{a\mathrm{d}\sigma}\right)^{p} P\left(\widehat{T}_{a\mathrm{d}\sigma}\right) \right] = \\ \frac{1}{2^{\frac{d'}{2}}} \int_{\Gamma} \quessianobb{s^{p}P(s)} \ \mathrm{d}\sigma(\xi) \ .
\end{multline*}
\end{lemma}
\end{description}

\begin{proof}{\bf (Proof of Lemma \ref{claimproductoconpolinomio})}\\
In order to ease the notation, we will write 
$$\cirquess := \frac{1}{2^{\frac{d'}{2}}} \int_{\Gamma} \quessianobb{\varphi} \ \mathrm{d}\sigma(\xi) \ .$$
and $S:=\topnorm$.
Let $\epsilon>0$ and $\{ g_{j} \}_{j\in \mathbb{N}}$ be a sequence of polynomials without constant terms that converges uniformly to the function $s^{p}$ in $[0,R]$.
For each $j \in \mathbb{N}$ we have
\begin{eqnarray}\label{descomposicionenlastresis}
& \ & \left\vert \left(\frac{\pi}{\alpha}\right)^{\frac{d}{2}} \mathrm{Tr}\left[ \left(S\right)^{p} P\left(S\right) \right] - \cirquess[s^{p}P(s)] \right\vert \ \leq \nonumber \\
& \ & \ \nonumber \\
& \ & \quad \left(\frac{\pi}{\alpha}\right)^{\frac{d}{2}} \bigg\vert \mathrm{Tr}\big[ \left(S^{p}-g_{j}(S)\right) P\left(S\right) \big] \bigg\vert + \left\vert \left(\frac{\pi}{\alpha}\right)^{\frac{d}{2}} \mathrm{Tr}\big[ g_{j}(S) P\left(S\right) \big] - \cirquess[g_{j}P] \right\vert \nonumber \\
& \ & \ \nonumber \\
& \ & \quad \qquad + \ \bigg\vert \cirquess[( g_{j}-s^{p} )P] \bigg\vert \ =: \mathring{I}_{1}+\mathring{I}_{2}+\mathring{I}_{3} \ . 
\end{eqnarray}
We first bound $\mathring{I}_{1}$. Let $\widetilde{P}$ be the polynomial whose coefficients are the absolute values of the coefficients of $P$. In this way $ |P(s)| \leq \widetilde{P}(s) $ for all $s \in [0,R]$ and hence $\big( \widetilde{P}-|P| \big)(S)$ is a positive operator, which in turn implies $\mathrm{Tr}\big( \left\vert P\left(S\right) \right\vert \big) \leq \mathrm{Tr}\big( \widetilde{P}\left(S\right) \big) $. As $\widetilde{P}$ is a polynomial without constant term, applying Proposition \ref{corolariopolinomio} we can find a constant $C$ such that $ \left(\frac{\pi}{\alpha}\right)^{\frac{d}{2}} \mathrm{Tr}\big(  \widetilde{P}\left(S\right) \big) < C $. We can now get
\begin{eqnarray}\label{desigualdadparaiunoenlaclaim}
\mathring{I}_{1} & = & \left(\frac{\pi}{\alpha}\right)^{\frac{d}{2}} \bigg\vert \mathrm{Tr}\big[ \left(S^{p}-g_{j}(S)\right) P\left(S\right) \big]  \bigg\vert \nonumber \\
& \ & \ \nonumber \\
& \leq & \left(\frac{\pi}{\alpha}\right)^{\frac{d}{2}} \norma{P\left(S\right)}{\scriptscriptstyle S_{1}( A^{2}_{\alpha}(\mathbb{B}_{n}) )} \norma{S^{p}-g_{j}(S)}{\scriptscriptstyle \mathcal{L}( A^{2}_{\alpha}(\mathbb{B}_{n}) )} \nonumber \\
& \ & \nonumber \\
& = & \left(\frac{\pi}{\alpha}\right)^{\frac{d}{2}} \mathrm{Tr}\big( \left\vert P\left(S\right) \right\vert \big) \norma{s^{p}-g_{j}}{L^{\infty}[0,R]} \nonumber \\
& \ & \ \nonumber \\
& \leq & \left(\frac{\pi}{\alpha}\right)^{\frac{d}{2}} \mathrm{Tr}\big( \widetilde{P}\left(S\right) \big) \norma{s^{p}-g_{j}}{L^{\infty}[0,R]} \ \leq \ C \norma{s^{p}-g_{j}}{L^{\infty}[0,R]} \ . \\ 
& \ & \ \nonumber
\end{eqnarray}

To estimate $\mathring{I}_{3}$ we first notice that
$$ \norma{\frac{ (g_{j}-s^{p} )P}{s^{p}}}{L^{\infty}[0,R]} \leq \norma{g_{j}-s^{p}}{L^{\infty}[0,R]} \norma{\frac{P}{s^{p}}}{L^{\infty}[0,R]} \ , $$
with $ \norma{\frac{P}{s^{p}}}{L^{\infty}[0,R]} $ bounded as $P$ has no constant term and $p \in (0,1)$.
We then apply inequality (\ref{desigualdadparausarenlaqconcirculo}) from {\it d)} in Proposition \ref{propoconloscuatroincisosdelaclaseC} to get
\begin{eqnarray}\label{desigualdadparaitresenlaclaim}
\mathring{I}_{3} & \leq &  \frac{\sigma(\mathrm{supp} (a))}{p^{\frac{d'}{2}}} \left( \sup\limits_{\xi \in \Gamma} \frac{a(\xi)}{1-|\xi|^{n+1}}  \right)^{p} \norma{\frac{ (g_{j}-s^{p} )P}{s^{p}}}{L^{\infty}[0,R]} \nonumber \\
& \ & \ \nonumber \\
&\leq & C \norma{g_{j}-s^{p}}{L^{\infty}[0,R]} \ .
\end{eqnarray}
Both (\ref{desigualdadparaiunoenlaclaim}) and (\ref{desigualdadparaitresenlaclaim}) are valid for all $j \in \mathbb{N}$, and the constants involved in each of them do not depend on $j$ or $\alpha$. Since $ \norma{g_{j}-s^{p}}{L^{\infty}[0,R]} $ tends to zero, we can choose $j^{\ast} \in \mathbb{N}$ large enough so that $\mathring{I}_{1}<\frac{\epsilon}{3}$ and $\mathring{I}_{3} < \frac{\epsilon}{3}$. Finally, for this fixed $j^{\ast}$, $g_{j^{\ast}}P$ is a polynomial without constant term, so by Proposition \ref{corolariopolinomio}
$$ \mathring{I}_{2} = \left\vert \left(\frac{\pi}{\alpha}\right)^{\frac{d}{2}} \mathrm{Tr}\big[ g_{j^{\ast}}(S) P\left(S\right) \big] - \cirquess[g_{j^{\ast}}P] \right\vert \ < \ \frac{\epsilon}{3} $$
for all $\alpha$ large enough, which together with (\ref{descomposicionenlastresis}) proves the claim.
\end{proof}

We now give a proof of the main result of the paper.

\begin{proof}[Proof of Theorem \ref{teoremaszego} $\lbrack$Szeg\"{o} limit theorem$\rbrack$]
Because $\varphi \in \mathcal{C}$, $\varphi(0)=0$ and we can find $p \in (0,1)$ such that $\frac{\varphi(s)}{s^{p}} \in C[0,R]$ and $\frac{\varphi(s)}{s^{p}}=0$ at $s=0$. Let $\{ p_{j} \}_{j\in \mathbb{N}}$ be a sequence of polynomials without constant terms that converges uniformly to $\frac{\varphi(s)}{s^{p}}$ in $[0,R]$. This in turn implies that $s^{p}p_{j}$ converges uniformly to $\varphi$ in $[0,R]$. We use once more the notation for $\mathring{Q}(\varphi)$ and $S$ that we defined in the proof of Lemma \ref{claimproductoconpolinomio} and follow a similar approach. Let $\epsilon>0$. For each $j \in \mathbb{N}$ we have
\begin{eqnarray}\label{descomposicionenlastresisenlademodeszego}
& \ & \left\vert \left(\frac{\pi}{\alpha}\right)^{\frac{d}{2}} \mathrm{Tr}\big( \varphi(S) \big) - \cirquess[\varphi] \right\vert \ \leq \nonumber \\
& \ & \ \nonumber \\
& \ & \qquad \left(\frac{\pi}{\alpha}\right)^{\frac{d}{2}} \bigg\vert \mathrm{Tr}\big( \varphi(S)-S^{p}p_{j}(S) \big)  \bigg\vert + \left\vert \left(\frac{\pi}{\alpha}\right)^{\frac{d}{2}} \mathrm{Tr}\big( S^{p}p_{j}(S) \big) - \cirquess[s^{p}p_{j}] \right\vert \nonumber \\
& \ & \ \nonumber \\
& \ & \qquad \qquad + \ \bigg\vert \cirquess[s^{p}p_{j}-\varphi] \bigg\vert \ =: \mathring{I}_{4}+\mathring{I}_{5}+\mathring{I}_{6} \ . 
\end{eqnarray}
For $\mathring{I}_{4}$:
\begin{eqnarray}\label{desigualdadparaicuatroenelteorema}
\mathring{I}_{4} & = & \left(\frac{\pi}{\alpha}\right)^{\frac{d}{2}} \bigg\vert \mathrm{Tr}\big( (\varphi-s^{p}p_{j})(S) \big)  \bigg\vert \nonumber \\
& \ & \ \nonumber \\
& = & \left(\frac{\pi}{\alpha}\right)^{\frac{d}{2}} \bigg\vert \mathrm{Tr}\left( S^{p} \big( \frac{\varphi}{s^{p}} -p_{j} \big)(S) \right) \bigg\vert \nonumber \\
& \ & \ \nonumber \\
& \leq & \left(\frac{\pi}{\alpha}\right)^{\frac{d}{2}} \norma{S^{p}}{\scriptscriptstyle S_{1}( A^{2}_{\alpha}(\mathbb{B}_{n}) )} \norma{\big( \frac{\varphi}{s^{p}} -p_{j} \big)(S)}{\scriptscriptstyle \mathcal{L}( A^{2}_{\alpha}(\mathbb{B}_{n}) )} \nonumber \\
& \ & \nonumber \\
& = & \left(\frac{\pi}{\alpha}\right)^{\frac{d}{2}} \mathrm{Tr}\big( S^{p} \big) \norma{\frac{\varphi}{s^{p}}-p_{j}}{L^{\infty}[0,R]} \ \leq \ C \norma{\frac{\varphi}{s^{p}}-p_{j}}{L^{\infty}[0,R]} \ ,
\end{eqnarray}
where the last inequality is obtained applying Lemma \ref{TRAZAPCHICA}.

For $\mathring{I}_{6}$, we employ inequality (\ref{desigualdadparausarenlaqconcirculo}) from {\it d)} in Proposition \ref{propoconloscuatroincisosdelaclaseC} to obtain
\begin{equation}\label{desigualdadparaiseisenelteorema}
\mathring{I}_{6} \ \leq \ C \norma{\frac{s^{p}p_{j}-\varphi}{s^{p}}}{L^{\infty}[0,R]} \ = \ C \norma{\frac{\varphi}{s^{p}}-p_{j}}{L^{\infty}[0,R]} \ .
\end{equation}
Both (\ref{desigualdadparaicuatroenelteorema}) and (\ref{desigualdadparaiseisenelteorema}) hold for any $j$ and are independent of $\alpha$, so we can pick $j^{\ast} \in \mathbb{N}$ such that $ \mathring{I}_{4} < \frac{\epsilon}{3} $ and $\mathring{I}_{6} < \frac{\epsilon}{3}$ for all $\alpha > 0$. Finally, for $j^{\ast}$ fixed, Lemma \ref{claimproductoconpolinomio} asserts $\mathring{I}_{5}<\frac{\epsilon}{3}$ for all $\alpha$ large enough. Applying all these inequalities in (\ref{descomposicionenlastresisenlademodeszego}) concludes the proof of the theorem.
\end{proof}

\section{Corollaries and follow up results}
\label{laseccionconloscorolarios}

We present 
some applications of the Szeg\"{o} limit theorem \ref{teoremaszego}, concretely, results describing assymptotically the distribution of the eigenvalues and the Schatten norms of $\topnorm$. Both of these results are equivalent to the ones obtained in \cite{SalvadorAlejandroUribe} and their proofs are analogues to the ones therein.\par

\subsection{Distribution of eigenvalues}
We begin by stating the equality in Theorem \ref{teoremaszego} in an alternate form. Addressing the definition of $\mathcal{Q}$, the right hand side of (\ref{eqteoremaszego}) equals
$$ \frac{1}{2^{\frac{d'}{2}}\Gamma(\frac{d'}{2})} \int_{\Gamma} \int_{0}^{\scriptscriptstyle \frac{a(\xi)}{(1-|\xi|^{2})^{n+1}}} \varphi(s) \left( \ln\frac{a(\xi)}{s(1-|\xi|^{2})^{n+1}} \right)^{\frac{d'}{2}-1} \ \frac{\mathrm{d}s}{s} \mathrm{d}\sigma(\xi) \ , $$
and interchanging the integrals order by Fubini's theorem, it can be rewritten as
\begin{multline*}
    \frac{1}{2^{\frac{d'}{2}}} \int_{0}^{\scriptscriptstyle \max\limits_{\xi \in \Gamma}  \frac{a(\xi)}{(1-|\xi|^{2})^{n+1}}} \varphi(s) \\
    \cdot \left[ \frac{1}{\Gamma(\frac{d'}{2})s} \int_{\scriptscriptstyle \left\{ \xi \in \Gamma \ : \ \frac{a(\xi)}{(1-|\xi|^{2})^{n+1}} \geq s \right\}} \left( \ln\frac{a(\xi)}{s(1-|\xi|^{2})^{n+1}} \right)^{\frac{d'}{2}-1} \ \mathrm{d}\sigma(\xi) \right] \mathrm{d}s \ .
\end{multline*}

If we define $\mathcal{D}_{a}(s)$ to be the function defined by what is inside the brackets in the last expression, then (\ref{eqteoremaszego}) can be restated as
\begin{equation}
\lim\limits_{\alpha \rightarrow \infty} \left( \frac{\pi}{\alpha} \right)^{\frac{d}{2}} \mathrm{Tr}\left( \varphi\left(\widehat{T}_{a\mathrm{d}\sigma}\right) \right) = \frac{1}{2^{\frac{d'}{2}}} \int_{0}^{\scriptscriptstyle \max\limits_{\xi \in \Gamma}  \frac{a(\xi)}{(1-|\xi|^{2})^{n+1}}} \varphi(s) \mathcal{D}_{a}(s) \ \mathrm{d}s \ .
\nonumber
\end{equation}

Notice that by doing this, we are now interpreting the asymptotics given in the Szeg\"{o} limit theorem as an integral with respect to a measure $\mathcal{D}_{a}(s) \ \mathrm{d}s $ that is absolutely continuous with respect to the Lebesgue measure in $\mathbb{R}$.\par

\begin{corollary}\label{corolariodedistribuciondeeigenvalores}
Let $I \subset (0, \max\limits_{\xi \in \Gamma}  \frac{a(\xi)}{(1-|\xi|^{2})^{n+1}} ]$ be a closed interval and denote by ${\mathcal{N}_{I}(\topnorm)}$ the number of eigenvalues of $\topnorm$ in $I$. Then
\begin{equation}
\lim\limits_{\alpha \rightarrow \infty} \left( \frac{\pi}{\alpha} \right)^{\frac{d}{2}} \mathcal{N}_{I}(\topnorm) = \frac{1}{2^{\frac{d'}{2}}} \int_{I} \mathcal{D}_{a}(s) \ \mathrm{d}s \ . \nonumber
\end{equation}
\end{corollary}

\subsection{Schatten norm estimates}
Let $p>0$. Recall that for $A $ a compact operator on $A^{2}_{\alpha}(\mathbb{B}_{n})$, $A$ is a member of the Schatten class $S_{p}(A^{2}_{\alpha}(\mathbb{B}_{n}))$ if and only if $|A|^{p}$ is of trace class. This is equivalent to having the quantity
$$ \norma{A}{S_{p}(A^{2}_{\alpha}(\mathbb{B}_{n}))} := \mathrm{Tr}((A^{\ast}A)^{\frac{p}{2}})^{\frac{1}{p}} $$
being finite. 
For $p\geq 1$, $ \norma{\cdot}{S_{p}( A^{2}_{\alpha}(\mathbb{B}_{n}) )} $ is a norm in $S_{p}(A^{2}_{\alpha}(\mathbb{B}_{n}))$, while when $0<p<1$ it is a quasi-norm.

\begin{corollary}
Let $\Gamma$ be isotropic or co-isotropic and $a \in C^{\infty}_{0}(\Gamma)$. For any $p>0$,
\begin{equation}
\lim\limits_{\alpha \rightarrow \infty} \left( \frac{\pi}{\alpha} \right)^{\frac{d}{2p}} \normaaa{ \topnorm }{S_{p}( A^{2}_{\alpha}(\mathbb{B}_{n}) )} = \left( \frac{1}{(2p)^{\frac{d'}{2}}} \int_{\Gamma} \left(\frac{|a(\xi)|}{(1-|\xi|^{2})^{n+1}} \right)^{p} \ \mathrm{d}\sigma(\xi) \right)^{\frac{1}{p}} \ . \nonumber
\end{equation}
\end{corollary}

\subsection{An example}
We illustrate some of the results given by presenting a concrete example.

For $r \in (0,1) $ let $\Gamma = r\mathbb{S}^{1} = \{ z \in \mathbb{C} : |z|=r \} $. If we use the parametrization $ \gamma: (0,1) \rightarrow \Gamma $, $\gamma(\theta) = re^{2\pi i \theta}$, then $\mathrm{d}\sigma = 2\pi \frac{r}{1-r^{2}} \mathrm{d}\theta$. Notice that in this case $d=1$ and $n=1$, thus $d'=1$. Let $a \equiv 1$ in $\Gamma$ and denote $ T_{r} := T_{1 \mathrm{d}\sigma } $ the associated Toeplitz operator. Then
\begin{eqnarray}
T_{r}f(z) & = & 2\pi r (1-r^{2})^{\alpha-1} (\alpha+1) \int_{0}^{1} K^{\alpha}(z,re^{2\pi i\theta})f(re^{2\pi i\theta}) \ \mathrm{d}\theta \ \nonumber
\end{eqnarray}
for any $ f \in A^{2}_{\alpha}(\mathbb{D}) $. 

Recall that $ e_{m}(z) = \delta_{m} z^{m} $ with $ \delta_{m} := \sqrt{\frac{\Gamma(m+\alpha+2)}{m! \ \Gamma(\alpha+2)}} $, $m \in \mathbb{N}_{0}$, is an orthonormal basis for $A_{\alpha}^{2}(\mathbb{D})$, and that the Bergman kernel can be expressed as an hypergeometric series via
$$ K^{\alpha}(z,w) = \frac{1}{(1-z\overline{w})^{\alpha+2}} = \sum_{j=0}^{\infty} \ \frac{\Gamma(\alpha+j+2)}{\Gamma(\alpha+2)\ j!} (z\overline{w})^{j} \ . $$
Using this we compute
\begin{eqnarray}
\frac{T_{r}e_{m}(z)}{\scriptstyle 2\pi r (1-r^{2})^{-\alpha+1} } & = & (\alpha+1) \delta_{m} \int_{0}^{1} \left( \ \sum_{j=0}^{\infty} \ \frac{\Gamma(\alpha+j+2)}{\Gamma(\alpha+2)\ j!} z^{j}r^{j}e^{-2\pi j i\theta} \right) r^{m} e^{2\pi m i\theta } \mathrm{d}\theta \nonumber \\
& & \nonumber \\
& = & (\alpha+1) \delta_{m} \sum_{j=0}^{\infty} \frac{\Gamma(\alpha+j+2)}{\Gamma(\alpha+2)\ j!} r^{m+j} z^{j} \int_{0}^{1} e^{2\pi(m-j)i\theta} \ \mathrm{d} \theta \nonumber \\
& & \nonumber \\
& = & (\alpha+1) r^{2m} \frac{\Gamma(\alpha+m+2)}{\Gamma(\alpha+2)\ m!} e_{m}(z) \ , \nonumber 
\end{eqnarray}

from where it follows that the eigenvalues of $\widehat{T}_{r} = \frac{1}{\sqrt{2\pi \alpha}} T_{r}$ are 
\begin{equation}\label{loseigenvaloresenelejemplo}
\lambda_{m} = \sqrt{\frac{2\pi}{\alpha}}(1-r^{2})^{\alpha-1} \frac{\Gamma(\alpha+m+2)}{\Gamma(\alpha+1)m!} r^{2m+1} \ , 
\end{equation}
$ m \in \mathbb{N}_{0}$.

\subsubsection{Operator norm}
We seek to analyze the behavior of the operator norm of $ \widehat{T}_{r} $ by estimating its largest eigenvalue. For any fixed $\alpha $, we notice that the sequence of quotients
$$ \frac{\lambda_{m}}{\lambda_{m-1}} = \frac{(\alpha+1)r^{2}}{m} + r^{2} $$
is monotonically decreasing to $r^{2}<1$. Thus, the largest eigenvalue of $\widehat{T}_{r}$ for this $\alpha$ is $\lambda_{m^{*}}$, where $m^{*}$ is the largest natural number such that $ \frac{(a+1)r^{2}}{m*} + r^{2} \geq 1 $. This is obtained when $ m^{*} = \left\lfloor (\alpha+1) \frac{r^{2}}{1-r^{2}} \right\rfloor $. Assume $\alpha$ is such that $ (\alpha+1) \frac{r^{2}}{1-r^{2}}  \in \mathbb{N}$. Then in terms of the Beta function
\begin{equation}
\frac{\Gamma(\alpha+m^{*}+2)}{\Gamma(\alpha+1)m^{*}!} = \mathrm{\mathbf{B}} \left( \alpha + 1 , (\alpha+1)\frac{r^{2}}{1-r^{2}} \right)^{-1} r^{-2} \ . \nonumber
\end{equation}
If $\alpha$ grows, so does $m^{*}=(\alpha+1) \frac{r^{2}}{1-r^{2}}$, and we can derive that
$$ \mathrm{\mathbf{B}} \left( \alpha + 1 , (\alpha+1)\frac{r^{2}}{1-r^{2}} \right) \sim \sqrt{\frac{2\pi}{\alpha + 1}}(1-r^{2})^{\alpha+1} r^{2m^{*}-1} $$
from the known asymptotic formula $ \mathrm{\mathbf{B}} \left( x,y \right) \sim \sqrt{2\pi} \frac{x^{x-\frac{1}{2}}y^{y-\frac{1}{2}}}{(x+y)^{x+y-\frac{1}{2}}}$ for large $x$ and $y$
. Substituting this in the expression for $\lambda_{m^{*}}$ in (\ref{loseigenvaloresenelejemplo}) we get
$$ \lambda_{m^{*}} \ \sim \ \sqrt{\frac{\alpha+1}{\alpha}} \frac{1}{(1-r^{2})^{2}} \ \sim \ \frac{1}{(1-r^{2})^{2}} $$
when $\alpha$ grows. This shows that the normalization proposed in (\ref{lanormalizacionprimeravez}) and Corollary \ref{corolariodelanormalizacion} is optimal in terms of the growth rate of $\alpha$, while also exhibiting the sharpness of the bound given in Remark \ref{remarkdespuesdelanormalizacion}.

\subsubsection{Eigenvalue distribution comparison}

Keeping with $\widehat{T}_{r}$ as defined above,  
given an interval $[t_{1},t_{2}] \subset \left(0, \frac{1}{(1-r^{2})^{2}}\right)$ Corollary \ref{corolariodedistribuciondeeigenvalores} implies
\begin{multline}\label{simplificacioncorolariodistribucion}
\lim\limits_{\alpha \rightarrow \infty} \sqrt{ \frac{\pi}{\alpha}} \mathcal{N}_{[t_{1},t_{2}]}\big( \widehat{T}_{r} \big) \ = \ \frac{\sqrt{8\pi} r}{1-r^{2}} \left[ \ \sqrt{ \ln\left( \frac{1}{(1-r^{2})^{n+1}t_{1}} \right) } \right. \\
- \left. \ \sqrt{ \ln\left( \frac{1}{(1-r^{2})^{n+1}t_{2}} \right) } \ \right] \ ,
\end{multline}
which is also written as an asymptotic equivalence in the form
\begin{equation}\label{simplificacioncorolariodistribucionasintota}
\mathcal{N}_{[t_{1},t_{2}]}\big( \widehat{T}_{r} \big) \ \sim \ \frac{\sqrt{8\alpha} r}{1-r^{2}} \left[ \ \sqrt{ \ln\left( \frac{1}{(1-r^{2})^{n+1}t_{1}} \right) } \ - \ \sqrt{ \ln\left( \frac{1}{(1-r^{2})^{n+1}t_{2}} \right) } \ \right] \ .
\end{equation}

For fixed values of $r \in (0,1)$ and $t_{1} < t_{2} \leq (1-r^{2})^{-2}$ we have numerically counted the number of eigenvalues of $\widehat{T}_{r}$ in $[t_{1},t_{2}]$ for growing values of $\alpha$ using (\ref{loseigenvaloresenelejemplo}). Two comparisons between these numerical outputs and the expected results given by (\ref{simplificacioncorolariodistribucion}) and (\ref{simplificacioncorolariodistribucionasintota}) are presented in the next two tables.
\begin{table}[H]
\small
\begin{center}
\begin{tabular}{|c|c|c|c|}
\hline
\quad $\alpha$ \quad & \quad $\mathcal{N}_{[\frac{16}{15},\frac{16}{9}]}\big( \widehat{T}_{\frac{1}{2}} \big)$ \quad & \quad R.H.S. of (\ref{simplificacioncorolariodistribucionasintota}) \quad & \quad $\sqrt{\frac{\pi}{\alpha}} \mathcal{N}_{[\frac{16}{15},\frac{16}{9}]}\big( \widehat{T}_{\frac{1}{2}} \big)$ \quad \\ \hline
100 & 14 & 13.47 & 2.4814 \\ \hline
500 & 30 & 30.13 & 2.378 \\ \hline
$10^{3}$ & 42 & 42.61 & 2.3541 \\ \hline
$5 \cdot 10^{3}$ & 95 & 95.29 & 2.3813 \\ \hline
$10^{4}$ & 134 & 134.76 & 2.3750 \\ \hline
$5 \cdot 10^{4}$ & 301 & 301.35 & 2.3859 \\ \hline
$10^{5}$ & 426 & 426.17 & 2.3877 \\ \hline
\end{tabular}
\caption{\small Comparison for $r=\frac{1}{2}$, $t_{1}=\frac{16}{15}$, $t_{2}=\frac{16}{9}$. The values in the second column were obtained by counting numerically
 for how many values of $m$ the eigenvalue $\lambda_{m}$ as given in (\ref{loseigenvaloresenelejemplo}) lies in $[\frac{16}{15},\frac{16}{9}]$.  As expected, the numbers on the rightmost column seem to be tending to the value of the R.H.S. of (\ref{simplificacioncorolariodistribucion}), approximately $ 2.3887 $ in this case. }
\end{center}
\end{table}

\begin{table}[H]
\small
\begin{center}
\begin{tabular}{|c|c|c|c|}
\hline
\quad $\alpha$ \quad & \quad $\mathcal{N}_{[\frac{2}{5},\frac{3}{5}]}\big( \widehat{T}_{\frac{1}{\sqrt{2}}} \big)$ \quad & \quad R.H.S. of (\ref{simplificacioncorolariodistribucionasintota}) \quad & \quad $\sqrt{\frac{\pi}{\alpha}} \mathcal{N}_{[\frac{2}{5},\frac{3}{5}]}\big( \widehat{T}_{\frac{1}{\sqrt{2}}} \big)$ \quad \\ \hline
100 & 5 & 5.60 & 0.8862 \\ \hline
500 & 12 & 12.52 & 0.9511 \\ \hline
$10^{3}$ & 18 & 17.71 & 1.0088 \\ \hline
$5 \cdot 10^{3}$ & 39 & 39.61 & 0.9775 \\ \hline
$10^{4}$ & 56 & 56.02 & 0.9925 \\ \hline
$5 \cdot 10^{4}$ & 125 & 125.28 & 0.9908 \\ \hline
$10^{5}$ & 177 & 177.17 & 0.9920 \\ \hline
\end{tabular}
\caption{\small Comparison for $r=\frac{1}{\sqrt{2}}$, $t_{1}=\frac{2}{5}$, $t_{2}=\frac{3}{5}$. The values in the second column were obtained by counting numerically
 for how many values of $m$ the eigenvalue $\lambda_{m}$ as given in (\ref{loseigenvaloresenelejemplo}) lies in $[\frac{2}{5},\frac{3}{5}]$. As expected, the numbers on the rightmost column seem to be tending to the value of the R.H.S. of (\ref{simplificacioncorolariodistribucion}), approximately $ 0.9930 $ in this case.}
\end{center}
\end{table}

\appendix

\section{Proof of Lemma \ref{TRAZAPCHICA}}\label{seccioncondemodelaclaimdepchica}

Throughout this section $a$ will be a positive symbol and $0<p<1$. Lemma \ref{TRAZAPCHICA}
is equivalent to there being some constant $C>0$ such that
\begin{equation}\label{equivalenciatrazapchica}
\mathrm{Tr}\big( T_{a\mathrm{d}\sigma}^{p} \big) \leq C \alpha^{p\left(n-\frac{d}{2} \right)+\frac{d}{2}} \ ,
\end{equation}
for all $\alpha$ large enough. We will prove this inequality.

Let $\widetilde{T}_{a\mathrm{d}\sigma}$ be the Berezin symbol of $T_{a\mathrm{d}\sigma}$, namely $$ \widetilde{T}_{a\mathrm{d}\sigma}(z) = \langle T_{a\mathrm{d}\sigma} k_{z} , k_{z} \rangle \ .$$
For any Hilbert space $H$ and any positive operator $A \in \mathcal{L}(H)$, the relationship $ \langle A^{p}x , x \rangle \leq \langle Ax , x \rangle^{p} $ is true for all $x$ unitary in $H$ whenever $0<p<1$ (Proposition 1.31 \cite{khezhuchiquito}). Combining this with the definition of $\widetilde{T}_{a\mathrm{d}\sigma}$ yields in
$$ \widetilde{T_{a\mathrm{d}\sigma}^{p}}(z) \leq \left(\widetilde{T}_{a\mathrm{d}\sigma}(z)\right)^{p} \ $$
for all $z \in \mathbb{B}_{n}$. Using the formula for the trace of an operator in terms of its Berezin symbol, we get

\begin{align}\label{desigualdaddelatrazadepchicausandoberezin}
    \mathrm{Tr}\left( T_{a\mathrm{d}\sigma}^{p} \right) & = c_{\alpha}\int_{\mathbb{B}_{n}} \widetilde{T_{a\mathrm{d}\sigma}^{p}}(z) \ \mathrm{d}\tau(z) \nonumber \\
    & \leq  c_{\alpha}\int_{\mathbb{B}_{n}} \left(\widetilde{T}_{a\mathrm{d}\sigma}(z)\right)^{p} \ \mathrm{d}\tau(z) \nonumber \\
    & = c_{\alpha}^{1+p}\int\limits_{\mathbb{B}_{n}} \left[ \int\limits_{\Gamma} \left( \frac{(1-|z|^{2})(1-|w|^{2})}{|1-\langle z,w \rangle |^{2}} \right)^{n+1+\alpha} \frac{a(w) \mathrm{d}\sigma(w)}{(1-|w|^{2})^{n+1}} \right]^{p} \ \mathrm{d}\tau(z)
\end{align}
where $\mathrm{d}\tau(z) = \frac{\mathrm{d}v(z)}{(1-|z|^{2})^{n+1}}$ is the invariant measure on $\mathbb{B}_{n}$.

\begin{myremark}\label{remarkmetricahiperbolica}
If $\varphi_{z}$ denotes the involutive automorphism in $\mathbb{B}_{n}$ that interchanges $z$ with the origin, then $$ \frac{(1-|z|^{2})(1-|w|^{2})}{|1-\langle z,w \rangle |^{2}} = 1-|\varphi_{z}(w)|^{2} \ .$$
Moreover $|\varphi_{z}(w)|$ is the pseudo-hyperbolic distance between $z$ and $w$, while the hyperbolic distance is $\beta(z,w) = \tanh^{-1}|\varphi_{z}(w)|$. This then implies
$$ \frac{(1-|z|^{2})(1-|w|^{2})}{|1-\langle z,w \rangle |^{2}} \leq 1 \ , $$
with equality if and only if $\beta(z,w)=0$, if and only if $z=w$ (see \cite{khezhu} Sections 1.2 and 1.5, particularly Lemma 1.2 and Proposition 1.21).
\end{myremark}

For $\epsilon > 0$ we now define 
$$ \mathcal{V}= \{ z \in \mathbb{B}_{n} \ : \ \exists w \in \Gamma \ , \ \beta(z,w)<\epsilon \} $$
a hyperbolic tubular neighborhood of $\Gamma$ in $\mathbb{B}_{n}$.

\begin{proposition}\label{lapropodelavecindadtubular}
The neighborhood $\mathcal{V}$ satisfies
\small
\begin{equation}\label{ecuacionenlaproposiciondelavecindadtubular}
c_{\alpha}^{1+p}\int\limits_{\mathbb{B}_{n}\setminus \mathcal{V}} \left[ \int\limits_{\Gamma} \left( \frac{(1-|z|^{2})(1-|w|^{2})}{|1-\langle z,w \rangle |^{2}} \right)^{n+1+\alpha} \frac{a(w) \mathrm{d}\sigma(w)}{(1-|w|^{2})^{n+1}} \right]^{p} \mathrm{d}\tau(z) = O\left( \alpha^{-\infty} \right) .
\end{equation}
\normalsize
\end{proposition}
\begin{proof}
Let $r(a)>0$ be such that $\beta(w,0)<r(a)$ for all $w \in \mathrm{supp}(a)$. There exists some constant $ C_{r(a)} >0$ that only depends on $r(a)$ such that $$ C_{r(a)}^{-1} \leq \frac{1-|u|^{2}}{1-|v|^{2}} \leq C_{r(a)}$$ for any $u,v \in \mathbb{B}_{n}$ with $\beta(u,v)<r(a)$ (Lemma 2.20 in \cite{khezhu}). Let $\rho<1$ be such that $ C_{r(a)}(1-\rho^{2})<1/2 $.
It suffices to prove that 
$$ \mathcal{J}_{1}:= \int\limits_{ (\mathbb{B}_{n}\setminus \mathcal{V}) \cap \{ z\in \mathbb{B}_{n} : |z|\geq \rho \} } \Big[ \ \ast \ \Big]^{p} \mathrm{d}\tau(z) \quad \mbox{and} \quad \mathcal{J}_{2}:= \int\limits_{ (\mathbb{B}_{n}\setminus \mathcal{V}) \cap \{ z\in \mathbb{B}_{n} : |z| < \rho \} } \Big[ \ \ast \ \Big]^{p} \mathrm{d}\tau(z) $$
are both $O\left( \alpha^{-\infty} \right)$, where $\ast$ represents what is inside the brackets in (\ref{ecuacionenlaproposiciondelavecindadtubular}).\par
If $z \in \mathbb{B}_{n}$, $w \in \mathrm{supp}(a)$, the invariance under automorphisms of the  Bergman metric leads to $ \beta(\varphi_{z}(w),z) = \beta(w,0) < r(a) $, and so
\begin{eqnarray}
\frac{(1-|z|^{2})(1-|w|^{2})}{|1-\langle z,w \rangle |^{2}} & = & \frac{1-|\varphi_{z}(w)|^{2}}{1-|z|^{2}} (1-|z|^{2}) \nonumber \\
& \ & \ \nonumber \\
& \leq & C_{r(a)}(1-|z|^{2}) \ . \nonumber
\end{eqnarray}
Using this we compute
\begin{eqnarray}
\mathcal{J}_{1} & \leq & C_{r(a)}^{(n+1+\alpha)p} \left[ \int\limits_{\Gamma} \frac{a(w) \mathrm{d}\sigma(w)}{(1-|w|^{2})^{n+1}}  \right]^{p} \int\limits_{ \{ z\in \mathbb{B}_{n} : |z| \geq\rho \} } (1-|z|^{2})^{\alpha p-(n+1)(1-p)} \mathrm{d}v(z) \nonumber \\
& \ & \ \nonumber \\
& \leq & C \left( C_{r(a)}(1-\rho^{2}) \right)^{\alpha p} \nonumber \\
& \ & \ \nonumber \\
& < & C \frac{1}{2^{\alpha p}} \ = \ O\left( \frac{1}{\alpha^{\infty}} \right) \nonumber
\end{eqnarray}
for all $\alpha$ large enough so that $ \alpha p-(n+1)(1-p) > 0 $.\par
If $ z \notin \mathcal{V} $ then $\beta(z,w) \geq \epsilon$ for all $w \in \Gamma$ and hence 
$$ 1-|\varphi_{z}(w)|^{2} \ = \ 1- \tanh^{2} \beta(z,w) \ \leq \ 1- \tanh^{2} \epsilon  \ ,$$
which we use to conclude that
\begin{eqnarray}
\mathcal{J}_{2} & \leq & ( 1- \tanh^{2} \epsilon )^{(n+1+\alpha)p} \left[ \int\limits_{\Gamma} \frac{a(w) \mathrm{d}\sigma(w)}{(1-|w|^{2})^{n+1}}  \right]^{p} \int\limits_{ \{ z\in \mathbb{B}_{n} : |z|<\rho \} }  \mathrm{d}\tau(z) \nonumber \\
& \ & \ \nonumber \\
& \leq & C(1-\tanh^{2}\epsilon)^{\alpha p} \ = \ O\left( \frac{1}{\alpha^{\infty}} \right) \ . \nonumber
\end{eqnarray}
\end{proof}

We take a moment to focus on the case $d=2n$, so $\Gamma $ is an open subset of $\mathbb{B}_{n}$ and $T_{a\mathrm{d}\sigma} = T_{a} $. Considering the result in last statement and using it together with the first inequality in (\ref{desigualdaddelatrazadepchicausandoberezin}) we can get
$$ \mathrm{Tr}\left( T_{a} \right) \leq c_{\alpha} \int_{\delta\mathbb{B}_{n}} \left( \widetilde{T_{a}}(z) \right)^{p} \ \mathrm{d}\tau(z) \ + \  O\left( \frac{1}{\alpha^{\infty}} \right) \ , $$
where $\delta<1$ is so that $\delta\mathbb{B}_{n} \supset \mathrm{supp}(a)$. In this case, the exponent of $\alpha$ in (\ref{equivalenciatrazapchica}) simplifies to $n$, which is the growth order of $c_{\alpha} = \frac{\alpha^{n}}{n!}(1+O(\alpha^{-1}))$. Therefore, it will suffice to show that $ \int_{\delta\mathbb{B}_{n}} \left( \widetilde{T_{a}}(z) \right)^{p} \ \mathrm{d}\tau(z) $ remains bounded when $\alpha $ grows. For this we recall that the Berezin transform is a contraction from $\mathcal{L}(A^{2}_{\alpha}(\mathbb{B}_{n}))$ to $L^{\infty}(\mathbb{B}_{n})$, so
\begin{eqnarray}
\int_{\delta\mathbb{B}_{n}} \left( \widetilde{T_{a}}(z) \right)^{p} \ \mathrm{d}\tau(z) & \leq & \normaaa{\widetilde{T_{a}}}{ L^{\infty}(\mathbb{B}_{n}) }^{p} \int_{\delta\mathbb{B}_{n}} \ \mathrm{d}\tau(z) \nonumber \\
& \ & \ \nonumber \\
& \leq & C \norma{T_{a}}{ \mathcal{L}(A^{2}_{\alpha}(\mathbb{B}_{n})) }^{p} \ . \nonumber
\end{eqnarray}
Since $d=2n$ implies $\Gamma$ being co-isotropic and $d'=0$, $\norma{T_{a}}{\mathcal{L}(A^{2}_{\alpha}(\mathbb{B}_{n}))} = \frac{\pi^{n}}{n!} \big\vert \big\vert \widehat{T}_{a}\big\vert \big\vert_{\scriptscriptstyle \mathcal{L}(A^{2}_{\alpha}(\mathbb{B}_{n}))} $ is bounded for all $\alpha>0$ by Corollary \ref{corolariodelanormalizacion}, thus proving the claim for $d=2n$.

Let $d\leq 2n-1$. From Proposition \ref{lapropodelavecindadtubular}, and rewriting the integral expression in (\ref{desigualdaddelatrazadepchicausandoberezin}) we conclude that proving our claim relies on bounding the integral
$$ \mathcal{I}_{p}(\alpha) := c_{\alpha}^{1+p}\int\limits_{\mathcal{V}} \left[ \int\limits_{\Gamma} e^{\alpha \ln (1-|\varphi_{z}(w)|^{2})} \frac{a(w) \ \mathrm{d}\sigma(w)}{|1-\langle z,w \rangle|^{2(n+1)}} \right]^{p} \ \frac{\mathrm{d}v(z)}{(1-|z|^{2})^{(n+1)(1-p)}} \ . $$
\ \\
We now give a parametrization for $\mathcal{V}$. By means of using a partition of unity subordinated to an open cover of the support of $a$ if necessary, we can assume without loss of generality that the support of $a$ is contained on the image of a single parametrization ${\gamma : B(0,1) \rightarrow \Gamma}$ of $\Gamma$. For each $p=\gamma(t) \in \Gamma $ let $\{ E_{1}(t), \sdots E_{2n-d}(t) \}$ be an orthonormal basis for the normal space ${ \mathrm{D}_{p}\iota\left( \mathrm{T}_{p}\Gamma \right)^{\perp} \subset \mathrm{T}_{p}\mathbb{B}_{n} }$ at $p$, so the ambient tangent space $\mathrm{T}_{p}\mathbb{B}_{n}$ splits as a direct sum 
$$\mathrm{T}_{p}\mathbb{B}_{n} =  \mathrm{D}_{p}\iota\left( \mathrm{T}_{p}\Gamma \right) \oplus \mathrm{span}\{ E_{1}(t), \sdots , E_{2n-d}(t) \}  \ .$$
Also let $\exp_{\gamma(t)}: \mathrm{T}_{p}\mathbb{B}_{n} \rightarrow \mathbb{B}_{n}$ denote the exponential map. As the Bergman ball is geodesically complete, $\exp_{\gamma(t)}$ has all of $\mathrm{T}_{p}\mathbb{B}_{n}$ as its domain. Denote by $B'(0,\epsilon) \subset \mathbb{R}^{2n-d}$ the ball of radius $\epsilon$ centered at the origin. We define
\begin{eqnarray}
\eta: B(0,1) \times B'(0,\epsilon) \subset \mathbb{R}^{d}\times \mathbb{R}^{2n-d} & \rightarrow & \mathbb{B}_{n} \nonumber \\
\eta(t,x) & = & \exp_{\gamma(t)}\left( \sum_{j=1}^{2n-d} x_{j}E_{j}(t) \right) \ . \nonumber
\end{eqnarray}

By choosing $\epsilon$ small enough, we can guarantee that:
\begin{itemize}
\item[i)] $\eta(t,x)$ is in a geodesic ball centered at $\gamma(t)$ and therefore the Bergman and radial distances between $\eta(t,x)$ and $\gamma(t)$ coincide, this is $$ \beta\left(\gamma(t), \eta(t,x)\right) = \left( \sum_{j=1}^{2n-d} x_{j}^{2} \right)^{\frac{1}{2}}  \ .$$
\item[ii)] The Bergman distance from $\eta(t,x)$ to $\Gamma$ is attained at $\gamma(t)$, this is $$ \beta\left(\gamma(t), \eta(t,x)\right) \leq \beta\left(\gamma(s), \eta(t,x)\right) $$ for all $\gamma(s) \in \Gamma$. Moreover $t$ is the only critical point of the function $ \beta\left(\gamma( \cdot ), \eta(t,x)\right) $ in $B(0,1)$, therefore the last inequality is strict whenever $s \neq t$. 
\end{itemize}

It follows that $\mathcal{V} = \eta( B(0,1) \times B'(0,\epsilon) )$ and
\begin{multline}\label{expresiontubularparametrizada}
\mathcal{I}_{p}(\alpha) := c_{\alpha}^{1+p}\int\limits_{B(0,1) \times B'(0,\epsilon)} \left[ \int\limits_{B(0,1)} e^{i \alpha \Psi((t,x),s)} \frac{a(\gamma(s)) \ \mathfrak{m}(s) \ \mathrm{d}s}{|1-\langle \eta(t,x),\gamma(s) \rangle|^{2(n+1)}} \right]^{p} \\ \cdot \frac{\mathfrak{n}(t,x) \ \mathrm{d}t \mathrm{d}x}{(1-|\eta(t,x)|^{2})^{(n+1)(1-p)}} \ ,
\end{multline}
where $ \mathfrak{m}(s) \mathrm{d}s = \mathrm{d}\sigma(w) $, $ \mathfrak{n}(t,x) \mathrm{d}t \mathrm{d}x = \mathrm{d}v(z) $, and we have introduced the phase function $ \Psi((t,x),s) = -i \ln \big(1-\big| \varphi_{\eta(t,x)}(\gamma(s))\big|^{2}\big)$. The observations made in {\rm ii)} and in remark \ref{remarkmetricahiperbolica} translate to that for each $(t,x) \in B(0,1) \times B'(0,\epsilon) $, the function $ \Psi((t,x),\cdot \ ) $ has its only critical point at $s=t$.

Applying the method of stationary phase (Theorem 7.7.5 in \cite{hormander}) yields in
\begin{multline}\label{faseestacionariapchica}
    \int\limits_{B(0,1)} e^{i \alpha \Psi((t,x),s)} \frac{a(\gamma(s)) \ \mathfrak{m}(s) \ \mathrm{d}s}{|1-\langle \eta(t,x),\gamma(s) \rangle|^{2(n+1)}} \\
    = \left( \frac{2\pi}{\alpha} \right)^{\frac{d}{2}} e^{i \alpha \Psi((t,x),t)} \frac{a(\gamma(t)) \ \mathfrak{m}(t)}{|1-\langle \eta(t,x),\gamma(t) \rangle|^{2(n+1)}} \cdot \frac{1+O\left(\frac{1}{\alpha}\right)}{\mathrm{det}(-i \mathrm{Hess}\Psi((t,x),t))^{\frac{1}{2}}} \ .
\end{multline}

The continuity of the functions involved and the compactness of the domain also imply that
$$ \frac{a(\gamma(t)) \ \mathfrak{m}(t) }{|1-\langle \eta(t,x),\gamma(t) \rangle|^{2(n+1)}} \cdot \frac{\mathrm{det}(-i \mathrm{Hess}\Psi((t,x),t))^{-\frac{1}{2}} \mathfrak{n}(t,x)}{(1-|\eta(t,x)|^{2})^{(n+1)(1-p)}} $$ is uniformly bounded for $(t,x) \in B(0,1) \times  B'(0,\epsilon)$. Implementing this and (\ref{faseestacionariapchica}) in the expression for $\mathcal{I}_{p}(\alpha)$ in (\ref{expresiontubularparametrizada}) we find that there exists some constant $C>0$ such that for all $\alpha$ large enough

\begin{eqnarray}
\frac{1}{c_{\alpha}^{1+p}} \ \mathcal{I}_{p}(\alpha) & \leq & C \left( \frac{2\pi}{\alpha} \right)^{\frac{dp}{2}} \int_{B'(0,\epsilon)} e^{i p \alpha \Psi((t,x),t)} \ \mathrm{d}x \nonumber \\
& \ & \ \nonumber \\
& = & C \left( \frac{2\pi}{\alpha} \right)^{\frac{dp}{2}} \int_{B'(0,\epsilon)} \Big(1-\tanh^{2}\beta\big(\gamma(t),\eta(t,x)\big) \Big)^{p\alpha} \ \mathrm{d}x  \nonumber \\
& \ & \ \nonumber \\
& = & C \left( \frac{2\pi}{\alpha} \right)^{\frac{dp}{2}} \int_{B'(0,\epsilon)} \Big(1-\tanh^{2}\left| x \right| \Big)^{p\alpha} \ \mathrm{d}x \ . \label{integralenterminosdelatangente}
\end{eqnarray}

We can find some constant $C>0$ such that $Cr^{2n-d-1} \leq (\tanh r)^{2n-d-1}$ for all $r \in (0,\epsilon)$. This can be done because $\frac{\tanh r}{r}$ tends to $1$ as $r$ approaches $0$. Integrating in polar coordinates:
\begin{eqnarray}
\int_{B'(0,\epsilon)} \Big(1-\tanh^{2}\left| x \right| \Big)^{p\alpha} \ \mathrm{d}x & = & (2n-d)\int_{0}^{\epsilon} r^{2n-d-1}(1-\tanh^{2}r)^{p\alpha} \ \mathrm{d}r \nonumber \\
& \leq & \frac{2n-d}{C} \int_{0}^{\epsilon} \tanh^{2n-d-1}r \ (1-\tanh^{2}r)^{p\alpha} \mathrm{d}r \ \nonumber \\
& = & \frac{n-\frac{d}{2}}{C} \int_{0}^{\tanh^{2} \epsilon} \rho^{n-\frac{d}{2}-1} \ (1-\rho)^{p\alpha -1} \ \mathrm{d}\rho
\nonumber \\
& \leq & \frac{n-\frac{d}{2}}{C} \ \mathrm{\mathbf{B}}\left( n-\frac{d}{2} ,p\alpha \right) 
\ , \label{cotaconfuncionbeta}
\end{eqnarray}
where we have used the change of variables $\rho=\tanh^{2}r$ and $\mathrm{\mathbf{B}}(z,w) $ is the Beta function. As a consequence of the Stirling's formula the Beta function has the asymptotic  approximation 
$$ \mathrm{\mathbf{B}}\left( n-\frac{d}{2} ,p\alpha \right) \sim \Gamma\left(n-\frac{d}{2}\right) (p\alpha)^{-n+\frac{d}{2}} $$
as $\alpha \rightarrow \infty$. Combining this and (\ref{cotaconfuncionbeta}) with (\ref{integralenterminosdelatangente}) and the fact that $ c_{\alpha} \sim \alpha^{n} $ we get
\begin{eqnarray}
\mathcal{I}_{p}(\alpha) & \leq & C \alpha^{n(1+p)-\frac{dp}{2}-n+\frac{d}{2}} \ = \  C \alpha^{p(n-\frac{d}{2})+\frac{d}{2}} \nonumber
\end{eqnarray}
which by Proposition \ref{lapropodelavecindadtubular} implies (\ref{equivalenciatrazapchica}), concluding the proof of Lemma \ref{TRAZAPCHICA}.

\section{The determinant of the Hessian}
\label{laseccionconloscalculosparaeldeterminante}

\subsection{The determinant as a block matrix}\label{subseccionochopuntouno}

Recall from (\ref{formamatricialdelquessianodelaphifactorizada}) that computing the determinant of the Hessian of $\varphi(t, \pmb{0})$ relies on computing that of $\mathcal{O}_{m-1}$, which we have defined as the $(m-1)d \times (m-1)d$ matrix that expressed as a $(m-1)\times (m-1)$ block matrix with blocks of size $d\times d$ has the tridiagonal form
\begin{equation}
\mathcal{O}_{m-1} = 
\begin{pmatrix}
2I & -I-i W & \scdots & 0 \\
-I+i W & 2I  & \scdots & 0 \\
\vdots & \vdots & \vdots & \vdots \\
0 & \scdots & 2I & -I-i W \\
0 & \scdots & -I+i W & 2I 
\end{pmatrix} . \ \ \nonumber
\end{equation}
Here $W:=G^{-1}H$, where $G$ and $H$ are the matrices encoding the Riemannian metric $g_{p}$ and the skew-symmetric form $h_{p}$ on $\mathrm{T}_{p}\Gamma$ in terms of a parametrization $\gamma$ of $\Gamma$ with $\gamma(t)=p$. This same block matrix appeared in \cite{SalvadorAlejandroUribe} in the Fock space situation.

Let $ \mathcal{R} $ stand for the ring of $d\times d$ complex matrices generated by $I$ and $W$. Then $\mathcal{O}_{m-1}$ can also be regarded as a $(m-1)\times (m-1)$ matrix with entries in $\mathcal{R}$. As $ \mathcal{R}$ is commutative, any such matrix possesses a well-defined determinant in $\mathcal{R}$, which we will denote $\mathrm{det}_{\mathcal{R}}$, and that possesses all the usual properties of determinants of matrices with complex entries.
The relation between this determinant and the usual one $\det \mathcal{O}_{m-1}$ is given in Theorem 1 in \cite{sylvester}. From that theorem it follows that
$$ \mathrm{det} \mathcal{O}_{m-1}= \mathrm{det}\big( \detr{m-1} \big) \ . $$
Notice that the right-hand side of the equation makes sense since $\mathrm{det}_{\mathcal{R}}\mathcal{O}_{m-1}$ is an element of $\mathcal{R}$, and thus a $d\times d$ matrix with complex entries. This relation was used in \cite{SalvadorAlejandroUribe} to compute $\det \mathcal{O}_{m-1}$. We cite the result therein expressing $ \detr{m-1} $ in terms of $W$, as all involved computations translate straightforwardly.

\begin{proposition}
Let $m \geq 2$. Then
\begin{equation}\label{laformapolinomialparaelvalordeldeterminanatedelaocaligraficaperomaschilo} \detr{m-1} = \sum_{\ell=0}^{ \lfloor \frac{m-1}{2} \rfloor } \binom{m}{2\ell+1}(-1)^{\ell}W^{2\ell} \ . 
\end{equation}
\end{proposition}

We interpret the geometric meaning of the matrix $W$. For every $p \in \Gamma$ define the endomorphism $ K_{p}: \mathrm{D}_{p}\iota ( \mathrm{T}_{p}\Gamma ) \rightarrow \mathrm{D}_{p}\iota ( \mathrm{T}_{p}\Gamma ) $ by
\begin{equation}\label{ladefiniciondelaksubp} 
K_{p} = \Pi_{p} \circ J \vert_{\mathrm{D}_{p}\iota ( \mathrm{T}_{p}\Gamma )} \ , 
\end{equation}
where $ \Pi_{p}: \mathrm{T}_{p}\mathbb{B}_{n} \rightarrow \mathrm{D}_{p}\iota ( \mathrm{T}_{p}\Gamma ) $ is the orthogonal projection and $ J \vert_{\mathrm{D}_{p}\iota ( \mathrm{T}_{p}\Gamma )} $ is the restriction of $J$ to $ \mathrm{D}_{p}\iota ( \mathrm{T}_{p}\Gamma ) $.

\begin{proposition}\label{laformadelamatrizw}
Consider $\gamma : B(0,1) \rightarrow \Gamma$ a parametrization of $\Gamma$ around $p$. Then $W$ is the matrix representation of the endomorphism $K_{p}$ with the basis for $\mathrm{D}_{p}\iota(\mathrm{T}_{p}\Gamma)$ given by $ \left\{ \mathrm{D}_{p}\iota \left( \tngvec{j}{t} \right) \right\}_{j=1,\dots , d} $.
\end{proposition}
\begin{proof}
The definition of $W$ is equivalent to $GW=H$, which in turn implies
\begin{eqnarray}\label{ecuaciondegporwenlademostracionsobrelaformadelamatrizw}
h_{jk} & = & \sum_{\ell=1}^{d} g_{j\ell}w_{\ell k} \nonumber \\
& = & \sum_{\ell=1}^{d} g_{p}\left( \tngvec{j}{t} , \tngvec{\ell}{t} \right) w_{\ell k} \nonumber \\
& = & \sum_{\ell=1}^{d} \mathring{b}_{p} \left( \mathrm{D}_{p}\iota\left( \tngvec{j}{t} \right) \ , \ \mathrm{D}_{p}\iota \left( \tngvec{\ell}{t} \right) \right) w_{\ell k} \nonumber \\
& = & \mathring{b}_{p} \left( \mathrm{D}_{p}\iota\left( \tngvec{j}{t} \right) \ , \ \sum_{\ell=1}^{d} w_{\ell k} \mathrm{D}_{p}\iota \left( \tngvec{\ell}{t} \right) \right)
\end{eqnarray}
for any $ 1\leq j,k \leq d $. On the other hand we know from the definition that $$ h_{jk} = \mathring{b}_{p} \left( \mathrm{D}_{p}\iota \left( \tngvec{j}{t} \right) \ , \ J \ \mathrm{D}_{p}\iota \left( \tngvec{k}{t} \right) \right) \ .$$
This together with (\ref{ecuaciondegporwenlademostracionsobrelaformadelamatrizw}) leads to
\begin{equation}
\mathring{b}_{p} \left( \mathrm{D}_{p}\iota \left( \tngvec{j}{t} \right) \ , \ J \ \mathrm{D}_{p}\iota \left( \tngvec{k}{t} \right) - \sum_{\ell=1}^{d} w_{\ell k} \mathrm{D}_{p}\iota \left( \tngvec{\ell}{t} \right) \right) = 0 \nonumber
\end{equation}
for any $ 1\leq j,k \leq d $. Since $ \left\{ \mathrm{D}_{p}\iota \left( \tngvec{j}{t} \right) \right\}_{j=1,\dots , d} $ is a basis for $\mathrm{D}_{p}\iota\left( \mathrm{T}_{p}\Gamma \right)$, we get 
$$ J \ \mathrm{D}_{p}\iota \left( \tngvec{k}{t} \right) - \sum_{\ell=1}^{d} w_{\ell k} \mathrm{D}_{p}\iota \left( \tngvec{\ell}{t} \right) \ \in \ \mathrm{D}_{p}\iota\left( \mathrm{T}_{p}\Gamma \right)^{\perp} \ , $$
and consequently
\begin{eqnarray}
\Pi_{p} \circ J \left( \mathrm{D}_{p}\iota \left( \tngvec{k}{t} \right) \right) & = & \sum_{\ell=1}^{d} w_{\ell k} \Pi_{p} \left( \mathrm{D}_{p}\iota \left( \tngvec{\ell}{t} \right) \right) \nonumber \\
& \ & \ \nonumber \\
& = & \sum_{\ell=1}^{d} w_{\ell k} \ \mathrm{D}_{p}\iota \left( \tngvec{\ell}{t} \right) \nonumber
\end{eqnarray}
for each $1\leq k \leq d$, which proves the statement.
\end{proof}

It can be verified from its definition in terms of $\Pi_{p}$ and $J$ that $K_{p}$ is skew-adjoint with respect to $\mathring{b}_{p}$. This implies that $W$ is skew-adjoint and hence all its nonzero eigenvalues are pure imaginary and come in conjugate pairs. Let $ \pm i \lambda_{k} $, $k=1,\sdots , r$ stand for said eigenvalues with possible repetitions due to multiplicity and satisfying
\begin{equation}\label{laenumeraciondelaslambdas}
0 < \lambda_{1} \leq \lambda_{2} \leq \scdots \leq \lambda_{r} \ .
\end{equation}

Let $P_{m-1}$ be the polynomial defined by the right-hand side of (\ref{laformapolinomialparaelvalordeldeterminanatedelaocaligraficaperomaschilo}). Then by the spectral mapping theorem the eigenvalues of $P_{m-1}(W)$ are 
\begin{eqnarray}
P_{m-1}(\pm i\lambda_{k}) & = & \sum_{\ell=0}^{ \lfloor \frac{m-1}{2} \rfloor } \binom{m}{2\ell+1}(-1)^{\ell}(\pm i\lambda_{k})^{2\ell} \ \nonumber \\
& \ & \ \nonumber \\
& = & \sum_{\ell=0}^{ \lfloor \frac{m-1}{2} \rfloor } \binom{m}{2\ell+1}\lambda_{k}^{2\ell} \ = \ \frac{(1+\lambda_{k})^{m}-(1-\lambda_{k})^{m}}{2\lambda_{k}} , \nonumber
\end{eqnarray} 
each with double the multiplicity as in $W$, together with the eigenvalue $P_{m-1}(0)= m$ with multiplicity $d-2r$ corresponding to the kernel of $W$. We use this to express $\det \mathcal{O}_{m-1} $ in terms of the eigenvalues of $W$.

\begin{proposition}\label{laproposiciondondeconcluimoseldeterminantedelaoenterminosdeloseigenvaloresdelaw}
Let $ \pm i \lambda_{k} $, $k=1,\sdots , r$ be the nonzero eigenvalues of W. Then
\begin{equation}\label{laecuacionconeldeterminantedelaoenterminosdeloseigenvaloresdelaw}
\sqrt{ \det \mathcal{O}_{m-1} } \ = \ m^{\frac{d}{2}-r} \prod_{k=1}^{r} \frac{(1+\lambda_{k})^{m}-(1-\lambda_{k})^{m}}{2\lambda_{k}} \ .
\end{equation}
\end{proposition}
\begin{proof}
\begin{align}
\det \mathcal{O}_{m-1} & = \det \big( \detr{m-1} \big) \nonumber \\
& = \det P_{m-1}(W) \nonumber \\
& = m^{d-2r} \prod_{k=1}^{r} \left( \frac{(1+\lambda_{k})^{m}-(1-\lambda_{k})^{m}}{2\lambda_{k}} \right)^{2} \ . \nonumber
\end{align}
\end{proof}

\subsection{The case of isotropic or co-isotropic manifolds}\label{subseccionochopuntodos}

In this section we describe the behavior of $K_{p}$ when $\Gamma$ is isotropic or co-isotropic, in order to express the eigenvalues $\pm i \lambda_{k}$ of $W$ needed to compute the determinant of $\mathcal{O}_{m-1}$.

Let $ (X,\omega) $ be a symplectic vector space of dimension $2n$, where the symplectic form is constructed using the formula $\omega( \cdot , \cdot ) = g( \cdot , J \cdot )$ with $J$ a complex structure on $X$ and $g$ an inner product on $X$ compatible with $J$. Let $V \subseteq X$ be a subspace of dimension $d$ and denote by $ V^{\omega} $ its symplectic complement as defined in (\ref{definiciondelcomplementosimplectico}) and by $V^{\perp}$ its orthogonal complement with respect to $g$.
These sets are also subspaces of $X$ and one can verify directly from the definitions that \begin{equation}\label{larelacionentreelcomplementosymplecticoyelortogonal}
V^{\perp} = J( V^{\omega} ) \qquad \mathrm{and} \qquad J(V^{\perp}) = V^{\omega} \ .
\end{equation}
Notice that this will imply $\mathrm{dim} V^{\omega} = \mathrm{dim} V^{\perp} = 2n-d$.

According to the definition of $K_{p}$, we define the endomorphism $K_{V}: V \rightarrow V$ by 
$$ K_{V} = \Pi_{V} \circ J \vert_{V} \ , $$
where $\Pi_{V}: X \rightarrow V$ is the orthogonal projection and $J \vert_{V}$ is the restriction of $J$ to $V$. We will describe these transformations, obtaining results that we will then interpret to the particular case of $K_{p}$. These calculations appear in \cite{SalvadorAlejandroUribe} with $\mathrm{T}_{p}\Gamma$  instead of $V$. 

\begin{proposition}\label{laproposicionsobreelkerneldelacomposiciondelajotaylaproyeccion}
$$\mathrm{ker} K_{V} = V^{\omega} \cap V \ .$$
\end{proposition}
\begin{proof}
Using the relation in (\ref{larelacionentreelcomplementosymplecticoyelortogonal}),
\begin{eqnarray}
\mathrm{ker} K_{V} & = & \ker \big( \Pi_{V} \circ J \vert_{V} \big) \nonumber \\
& = & J^{-1}\big( \mathrm{ker} \ \Pi_{V} \big) \cap V \ \nonumber \\
& = & J^{-1} \big( V^{\perp} \big) \cap V \nonumber \\
& = & V^{\omega} \cap V \ . \nonumber
\end{eqnarray}
\end{proof}

Recall from Definition \ref{definiciondeisotropicoyco-isotropico} that $V$ is isotropic if $V\subseteq V^{\omega}$ and co-isotropic if $V\supseteq V^{\omega}$. The next corollary follows straightforwardly.

\begin{corollary}
If $V$ is isotropic, then $K_{V}$ is identically zero.
\end{corollary}

This fully describes the isotropic case. We turn our attention to the co-isotropic case.

\begin{proposition}
If $V$ is co-isotropic, then the image of $K_{V}$ satisfies
$$ K_{V} \big( V \big) = J(V) \cap V \ . $$
\end{proposition}
\begin{proof}
Assume $\xi = K_{V} v$ for some $v \in V$. Then $ Jv = \xi + u $ for some $u \in V^{\perp} = J \big( V^{\omega} \big)$. Next we can find some $w \in V^{\omega}$ such that $u=Jw$ and hence $\xi = J(v-w) \in J(V)$, since $ w \in V^{\omega} \subseteq V $ because $V$ is co-isotropic.\par
For the reverse inclusion, if $\xi \in J(V)\cap V$, then $\xi=Jv$ for some $v \in V$, and therefore ${ \xi = \Pi_{V}\xi = \Pi_{v} Jv  = K_{V}v}$.
\end{proof}

\begin{corollary}
If $V$ is co-isotropic then it splits as a direct sum 
\begin{equation}\label{elsplittingenelcasoco-isotropico}
V = K_{V} \big( V \big) \oplus V^{\omega} \ .
\end{equation}
\end{corollary}
\begin{proof}
By Proposition \ref{laproposicionsobreelkerneldelacomposiciondelajotaylaproyeccion} $\mathrm{ker} K_{V} = V^{\omega} $, and thus
$$ \dim \big( \ker K_{V} \big) = \dim V^{\omega} = 2n-d \ . $$
On the other hand, by the dimension theorem
we have $\dim \big( \ker K_{V} \big) = d- \dim K_{V} \big( V \big) $. Comparing both expressions yields $\dim K_{V} \big( V \big) =2d-2n $. It follows that
$$ \dim V^{\omega} = 2n-d \qquad \mathrm{and} \qquad \dim K_{V} \big( V \big) = 2d-2n \ ,$$
which together add up to $d=\dim V$. Finally, by the previous proposition $ K_{V} \big( V \big) \subseteq J(V) $, and therefore
\begin{equation}
K_{V} \big( V \big) \cap V^{\omega} \ \subseteq \ J(V) \cap J(V^{\perp}) \ = \ \emptyset \ . \nonumber
\end{equation}
\end{proof}

\begin{remark}
We have so far shown that
\begin{equation}\label{laecuacionenlaremarksobrelaparidaddelrangodelak}
 \dim K_{V} \big( V \big) = \left\{ \begin{array}{cl} 0 & V \mbox{ isotropic} \\
2(d-n) & V \mbox{ co-isotropic} \end{array} \right. \ . \nonumber
\end{equation}
In general, the rank of $K_{V}$ is always pair regardless of $V$ being isotropic or co-isotropic because $K_{V}$ is skew-adjoint. 
In order to maintain consistency with the definition we gave for the eigenvalues of $K_{p}$ in (\ref{laenumeraciondelaslambdas}), we will use $r$ to denote half the rank of $K_{V}$.
\end{remark}

Keeping up with $V$ co-isotropic, consider the action of $K_{V}$ on the splitting (\ref{elsplittingenelcasoco-isotropico}) of $V$: If $\xi \in V^{\omega}$, then $J\xi \in JV^{\omega} = V^{\perp}$ and hence $ K_{V}\xi = \Pi_{V} J \xi = 0 $. Alternatively, if $\xi \in K_{V}\big( V \big) \subseteq JV$, then $J\xi \in V$, thus in this scenario $ K_{V}\xi = \Pi_{V}J\xi = J\xi $. From this we conclude that in the co-isotropic case, $K_{V}$ behaves like $J$ on $K_{V}\big( V \big)$ and like the zero function on $V^{\omega}$. We can therefore find bases for $K_{V}\big( V \big)$ and $V^{\omega}$ such that $K_{V}$ has the matrix representation
$$ \begin{pmatrix}
K_{V}
\end{pmatrix} = \begin{pmatrix}
0 & -I_{r\times r} & 0 \\
\phantom{-}I_{r \times r} & 0 & 0 \\
0 & 0 & \ \ 0 \ \ 
\end{pmatrix} \ , $$ 

where $I_{r \times r}$ is the identity matrix of size $r$. To obtain the eigenvalues of $K_{V}$, we now compute the zeros of

\begin{eqnarray}
\det \left[ \begin{pmatrix} K_{V} \end{pmatrix} - \lambda I_{d \times d} \right] & = & \det \begin{pmatrix}
-\lambda I_{r \times r} & -I_{r\times r} & 0 \\
I_{r\times r} & -\lambda I_{r\times r} & 0 \\
0 & 0 & -\lambda I_{(d-2r) \times (d-2r)}
\end{pmatrix} \nonumber \\
& \ & \ \nonumber \\
& = & \det \begin{pmatrix}
-\lambda I_{r \times r} & -I_{r\times r} \\
I_{r\times r} & -\lambda I_{r\times r}
\end{pmatrix} \det \begin{pmatrix}
-\lambda I_{(d-2r) \times (d-2r)}
\end{pmatrix} \nonumber \\
& \ & \ \nonumber \\
& = & (\lambda^{2}+1)^{r} (-\lambda)^{d-2r} \ . \nonumber
\end{eqnarray}

It follows that the eigenvalues of $K_{V}$ are $\pm i$ each with multiplicity $r$ and the zero eigenvalue with multiplicity $d-2r$.

We return to the special case of the symplectic space $ \big( \mathrm{T}_{p}\mathbb{B}_{n} , \omega_{p} \big) $, with $ \mathrm{D}_{p}\iota ( \mathrm{T}_{p}\Gamma ) \subset \mathrm{T}_{p}\mathbb{B}_{n} $ and $K_{p}$ as defined in (\ref{ladefiniciondelaksubp}). The results exposed in this subsection can be summarized in the next statement.

\begin{proposition}\label{laproposicionqueresumelasecciondegeometriasymplectica}
If $\Gamma$ is isotropic, then $K_{p}$ has no nonzero eigenvalues, this is, $r=0$. If $\Gamma$ is co-istropic, then the eigenvalues of $K_{p}$ are $\pm i$ each with multiplicity $r= d-n $ and the zero eigenvalue with multiplicity $ 2n-d $.
\end{proposition}

We conclude the section with the result we use in Section \ref{laseccionconelcomportamientoasintotico} expressing the value of $ \det \mathcal{O}_{m-1} $. It follows from (\ref{laecuacionconeldeterminantedelaoenterminosdeloseigenvaloresdelaw}) in Proposition \ref{laproposiciondondeconcluimoseldeterminantedelaoenterminosdeloseigenvaloresdelaw}, when taking into account the form of the eigenvalues $i\lambda_{k}$ described in the last proposition.

\begin{proposition}\label{elcorolariodondeconlcuimoseldeterminantedelaoexplicito}
Let $\Gamma$ be isotropic or co-isotropic. Then
\begin{equation}
\sqrt{ \det \mathcal{O}_{m-1} } \ = \ \left\{ \begin{array}{cl}
m^{\frac{d}{2}} & \Gamma \mbox{ isotropic} \\
2^{(d-n)(m-1)}m^{n-\frac{d}{2}} & \Gamma \mbox{ co-isotropic}
\end{array} \right. \ . \nonumber
\end{equation}
\end{proposition}

\end{document}